\documentclass[11pt,a4paper]{article}

\usepackage[english]{babel}
\usepackage[utf8]{inputenc}
\usepackage[T1]{fontenc}
\usepackage{textcomp}

\usepackage{indentfirst,latexsym,bm}
\usepackage{amsmath,amsfonts,amssymb,mathrsfs}
\usepackage{mathtools}
\mathtoolsset{showonlyrefs}
\usepackage{pifont}
\usepackage{amsthm}
\usepackage[bottom,flushmargin]{footmisc}
\allowdisplaybreaks[4]

\usepackage{array,multirow,graphicx}
\usepackage{picinpar}
\usepackage{setspace}
\usepackage{booktabs,threeparttable}
\usepackage[shortlabels]{enumitem}
\usepackage{authblk}
\usepackage[top=2cm,bottom=2cm,outer=2cm,inner=2cm]{geometry}
\usepackage[labelfont={bf,small},textfont={small}]{caption}
\usepackage{subcaption}
\usepackage{comment,appendix,multicol,tikz}
\usepackage{microtype}

\usepackage[textsize=footnotesize]{todonotes}
\usepackage[numbers]{natbib}
\usepackage{xcolor}
\usepackage[colorlinks=true,citecolor=blue,urlcolor=blue,linkcolor=blue]{hyperref}

\makeatletter
\def\namedlabel#1#2{%
  \begingroup
  #2%
  \def\@currentlabel{#2}%
  \phantomsection\label{#1}%
  \endgroup
}
\makeatother

\numberwithin{equation}{section}
\numberwithin{figure}{section}

\newcommand{\email}[1]{\href{mailto:#1}{\nolinkurl{#1}}}

\theoremstyle{plain}
\newtheorem{theorem}{Theorem}[section]
\newtheorem{lemma}[theorem]{Lemma}
\newtheorem{corollary}[theorem]{Corollary}
\newtheorem{proposition}[theorem]{Proposition}

\theoremstyle{definition}

\newtheorem{remark}[theorem]{Remark}

\newtheorem{assumptionH}{Assumption}

\numberwithin{equation}{section}

\newcommand{\E}{\mathbb{E}}

\newcommand{\pr}{\mathbb{P}}

\newcommand{\D}{\mathbf{D}}

\newcommand{\ang}[1]{\langle #1\rangle}

\begin{document}

\title{\bf{Two-sided estimates for two-point local time tails of branching random walks: a Green function criterion}}

\author{Shuwen Lai}
\date{}

\begingroup
\renewcommand{\thefootnote}{}
\footnotetext{%
\textit{Address:} School of Statistics and Data Science, Nankai University, Tianjin 300071, China.\\
\textit{Email:} laishuwen97@gmail.com
}
\endgroup

\maketitle

\begin{abstract}
This paper is devoted to the study of the joint tails of local times at two distinct sites for branching random walks (BRWs).
Our principle results establish a general Green-function criterion for two-sided exponential tail bounds with a common spatial factor given by a three-point Green function convolution, and further identify a strict inequality that must be satisfied by the corresponding exponential constants. 
In particular, we verify this criterion for subcritical BRWs on finitely generated groups, for critical BRWs on non-elementary hyperbolic groups, and for critical BRWs on $\mathbb{Z}^d$ with $d\geq 5$.
As applications, we revisit the growth rate of visited sites on spheres in hyperbolic groups obtained by Sidoravicius et al. [Commun. Pure Appl. Math. 2023], and also establish lower bounds for the probabilities of visiting finite sets in Euclidean lattices. 
\end{abstract}

\section{Introduction and main results}
A branching random walk (BRW) is a particle system with genealogy described by a Galton--Watson tree and spatial positions obtained by assigning independent increments along the genealogical edges.
As a basic quantify describing the spatial occupation of a BRW, the local time at a site $x$, denoted by $L(x)$, counts the total number of particles that ever visit $x$.
Motivated by the close connection to various fields, it has been has been widely studied in different settings; e.g. spatial epidemics models, super-Brownian motion, and the Abelian sandpile model; see \cite{LZ10, AHJ21}.
In this paper, we focus on the tail probabilities of local times of BRWs.

To the best of our acknowledge, the existing literature mainly focuses on limit theorems and scaling properties of local times and related occupation statistics; see, e.g. \cite{BZ07, LZ10, LZ11, LGL15}.
Recently, their distribution properties have attracted a large number of attention.
For instance, the tail behaviour of one-point local times was studied by Angel et al. \cite{AHJ21} for critical BRWs on $\mathbb{Z}^d$.
Gao and Lai \cite{GL26} subsequently established two-sided exponential estimates for one-point local time tails of subcritical BRWs on finitely generated groups and critical BRWs on non-elementary hyperbolic groups, with the spatial dependence characterized by the Green function of the underlying random walk (RW).
For further background on local times and related occupation statistics, we refer to \cite{GL26} and the references therein.

The one-point estimates can only describe the marginal distributions of local times, but do not account for the cross-dependence among visits to distinct sites. 
This leads naturally to the joint tail
\begin{equation} \label{eqn.01}
     \mathbb{P}_a \bigl(L(x)\geq n, L(y)\geq m\bigr), \quad n,m \geq1.
\end{equation}
As we illustrate later, such joint distributions would play an important role in second-moment arguments for the range, trace, and visiting probabilities of BRWs.
For a critical BRW on $\mathbb{Z}^d$ with $d\geq 5$, an upper bound on the probability of visiting two disjoint finite sets was obtained in \cite{Zhu17}, and recently, an upper bound on the joint tail probability of occupation times in finitely many balls was established in \cite{ASS26} in terms of branching capacity.
In low dimensions, Chen and Lai \cite{CL26} obtained two-sided estimates for critical BRWs on $\mathbb{Z}^d$ with $d\leq 3$.

In this paper, we aim to establish a general Green function criterion that yields two-sided exponential estimates for \eqref{eqn.01}, with the same spatial factor in the upper and lower bounds, uniformly in the starting point $a$, the target sites $(x,y)$, and the local time levels $(n,m)$.
This criterion provides a unified mechanism for obtaining these estimates in different geometric settings, including finitely generated groups and non-elementary hyperbolic groups. 
In contrast to the solely upper bound in the existing literature, our results not only provide matching lower bounds, but also identify the spatial factor governing the joint tail.

\subsection{Main results}
Let $\Gamma$ be a finitely generated group  with identity $e$, and let $\mu$ be an admissible symmetric probability measure on $\Gamma$, where \textsl{admissibility} means that its support generates $\Gamma$ as a semigroup.
Consider a BRW on $\Gamma$ started at $a\in \Gamma$ with step distribution $\mu$ and offspring distribution $\nu$ on $\mathbb{N}_0$.
Write $p_n(x,y)$ for the $n$-step transition probabilities of the underlying RW, and define the Green function by
\begin{equation*}
	G_s(x,y):=\sum_{n=0}^{\infty}s^n p_n(x,y),\quad s\geq 0.
\end{equation*}
Since $\mu$ is admissible, the series $G_s(x,y)$ either converge or diverge simultaneously for all $x, y\in \Gamma$; see \cite[Lemma~1.7]{Woe00}.
Let $R\geq 1$ be the common radius of convergence of the Green functions and let
\begin{equation*}
      r:=\sum_{k\geq0}k\nu(k)>0.
\end{equation*}
We refer to $r<R$ as the subcritical case and $r=R$ as the critical case.
Throughout the paper, we impose the following assumption on the offspring distribution $\nu$.

\begin{assumptionH} \label{ass:offspring}
	$\nu(0)+\nu(1)<1$ and there exists some $c>0$ such that $\sum_{k=0}^{\infty} e^{ck}\nu(k)<\infty$.
\end{assumptionH}

Assumption~\ref{ass:offspring} ensures that genuine branching occurs with positive probability, and that the offspring distribution $\nu$ has a finite exponential moment and hence finite moments of all orders.

For $s>0$ such that $G_s<\infty$, and for $a,x,y\in \Gamma$, define 
\begin{equation*}
     \Phi_s(a;x,y):= \sum_{z\in\Gamma} G_s(a,z)G_s(z,x)G_s(z,y).
\end{equation*}
We refer to this convolution-type quantity as the three-point Green function convolution.
Here $z$ represents a possible site at which a particle branches, with its descendants potentially visiting the target sites $x$ and $y$.

\begin{theorem} \label{thm:two-point-tail-estimate}
Under Assumption~\ref{ass:offspring}, suppose that $G_r(x,y)<\infty$ for all $x,y\in \Gamma$, and that there exists a constant $C>0$ such that, for all $a,x,y \in \Gamma$,
\begin{gather}
    \sum_{z\in\Gamma}G_r(a,z)^2 G_r(z,x) \leq C G_r(a,x), \tag{G1} \label{eq:criterion-triangle}\\
    \sum_{z\in\Gamma}G_r(a,z)G_r(z,x)\Phi_r(z;x,y)  \leq C\Phi_r(a;x,y). \tag{G2}  \label{eq:criterion-mixed}
\end{gather}
Then, $\Phi_r(a;x,y)<\infty$ for all $a,x,y\in \Gamma$, and there exist constants $c_1,c_2,C_1,C_2>0$ such that, for all $a,x,y\in\Gamma$ with $x\neq y$ and all integers $n,m \geq 1$,
\begin{equation}  \label{eq:two-point-tail-estimate}
    c_1e^{-c_2(n+m)}\Phi_r(a;x,y) \leq \pr_a\left(L(x)\geq n,L(y)\geq m\right)  \leq C_1e^{-C_2(n+m)} \Phi_r(a;x,y).
\end{equation}
Moreover, any such choice of constants necessarily satisfies $c_2>C_2$.
\end{theorem}

Theorem \ref{thm:two-point-tail-estimate} shows that the same spatial factor $\Phi_r(a;x,y)$ appears in both the upper and lower bounds, while the dependence on the local time levels $n$ and $m$ is described by exponential factors in $n+m$.
The quantity $\Phi_r(a;x,y)$ reflects the contribution of common ancestors that have descendants visiting $x$ and $y$. 
In particular, the probability of visiting both target sites is comparable to $\Phi_r(a;x,y)$, and the conditional joint tail probabilities, given that both sites are visited, admit a uniform two-sided exponential bound in $n+m$.
A direct application of the Cauchy--Schwarz inequality to the one-point tail estimates in \cite{GL26} yields an upper bound with spatial factor $(G_r(a,x)G_r(a,y))^{1/2}$, whereas $\Phi_r(a;x,y)$ captures the additional cross-dependence between visits to distinct sites.

Theorem~\ref{thm:two-point-tail-estimate} reduces the problem of establishing two-sided exponential bounds for the joint tail of $L(x)$ and $L(y)$ to verifying the conditions \eqref{eq:criterion-triangle} and \eqref{eq:criterion-mixed} on the Green function.
We next verify these conditions in three different settings. 
First, all conditions are readily verified for a subcritical BRW on any finitely generated group $\Gamma$.

\begin{corollary} \label{cor:subcritical-two-point-tail-estimate}
When $r\in(0,R)$, the two-sided bound \eqref{eq:two-point-tail-estimate} holds for any $\Gamma$, and necessarily $c_2>C_2$.
\end{corollary}

In the critical regime ($r=R$), the verification of \eqref{eq:criterion-triangle} and \eqref{eq:criterion-mixed} requires additional geometric input.
When $\Gamma$ is a non-elementary hyperbolic group endowed with a word metric $d_{\Gamma}$, we say that the probability measure $\mu$ on $\Gamma$ has \textsl{superexponential tails} if, for every $\lambda\geq 0$,
\begin{equation*}
	\sum_{g\in\Gamma}e^{\lambda d_{\Gamma}(e,g)}\mu(g)<\infty.
\end{equation*} 
The assumption of superexponential tail ensures Ancona inequalities at $R$, which allow us to verify \eqref{eq:criterion-triangle} and \eqref{eq:criterion-mixed}.
We write $f\asymp g$ if  $cg\leq f\leq Cg$ for some constants $c,C>0$.

\begin{corollary}  \label{cor:critical-two-point-tail-estimate-hyperbolic}
Assume that $\Gamma$ is a non-elementary hyperbolic group and $\mu$ has superexponential tails.
If $r=R$, then the two-sided bound \eqref{eq:two-point-tail-estimate} holds, and necessarily $c_2>C_2$.
Moreover,
\begin{equation*}
	\Phi_R(a;x,y) \asymp G_R(a, \pi(a)) G_R(\pi(a), x) G_R(\pi(a), y),
\end{equation*}
where $\pi(a)$ is a nearest point projection of $a$ onto a geodesic segment $[x,y]$.
\end{corollary}

We finally consider the critical BRW on Euclidean lattices. 
In this setting, we always have $R=1$.
In particular, when $d\geq 5$, the polynomial decay of the Green function given in \cite[Theorem~4.3.1]{LL10} is sufficiently fast to ensure that both the conditions \eqref{eq:criterion-triangle} and \eqref{eq:criterion-mixed} hold.
For $z\in\mathbb{Z}^d$, write $\ang{z}:=1\vee |z|$, where $|\cdot|$ is the $\ell^1$-norm.

\begin{corollary} \label{cor:critical-euclidean-two-point-tail-high}
Assume that $\Gamma=\mathbb{Z}^d$ with $d\geq 5$ and the underlying RW is simple.
If $r=R=1$, then the two-sided bound \eqref{eq:two-point-tail-estimate} holds, and necessarily $c_2>C_2$.
Moreover,
 \begin{equation*}
	\Phi_1(a;x,y) \asymp l_1^{-d+4} l_2^{-d+2},
\end{equation*}
where $l_1\leq l_2\leq l_3$ is the increasing rearrangement of $\ang{x-a}$, $ \ang{y-a}$ and $\ang{x-y}$.
\end{corollary}

In dimensions $d\leq 3$, the two-point local time tails decay polynomially in the local time levels, and their spatial dependence does not involve the mutual distance between the two target sites; see \cite{CL26}.
By contrast, for $d\geq 5$, the tails exhibit exponential decay in the local time levels, with a separate spatial factor that depends on the mutual distance.
The case $d=4$ is not covered here. We will consider it in \cite{Lai26} as a critical case.

As applications, we revisit the growth rate of visited sites on spheres in hyperbolic groups established in \cite{SWX23} (see Corollary~\ref{cor:critical-visited-sphere-growth}), and establish a general lower bound for visiting finite sets in $\mathbb{Z}^d$ (see Proposition~\ref{prop:critical-euclidean-finite-set-visiting}). 
For balls contained in coordinate subspaces, this lower bound revisits the corresponding lower bounds in \cite{Zhu26} under our offspring assumptions (see Corollary~\ref{cor:critical-euclidean-flat-ball-visiting}).

\subsection{Outline of proofs}
Our proof adapts the diagrammatic approach of \cite{AHJ21, GL26}.
For the upper bound, we expand mixed local time moments into diagrammatic sums and derive recursive bounds on their diagrams.
Decomposing a diagram yields three types of contributions: pairs of one-point diagrams, mixed pairs of one-point and  two-point diagrams, and pairs of two-point diagrams. 
The latter two cases introduce convolutions involving $G_r\Phi_r$ and $\Phi_r^2$, respectively.
The main task is to bound these terms while retaining the spatial factor $\Phi_r$.

For the lower bound, we extract skeletons with two labelled chains descending from a common branching vertex.
Their contribution yields mixed factorial moment lower bounds with the same spatial factor $\Phi_r$.
Applying the Paley--Zygmund inequality to $\binom{L(x)}{n}\binom{L(y)}{m}$, together with the mixed moment upper bounds, we obtain  the joint lower tail estimate.

Using the branching property, we further show that the exponential constants in any corresponding two-sided estimate must necessarily satisfy $c_2>C_2$.

\smallskip
{\it\textbf{Organization of this paper.}} 
Section~\ref{sec:preliminaries} introduces notation and basic facts. 
Section~\ref{sec:diagrammatic-expansion} develops the diagrammatic expansions and derive recursive bounds. 
Section~\ref{sec:general-two-point-criterion} establishes a general criterion for two-point tails and proves the subcritical case.
Sections~\ref{sec:hyperbolic-case} and \ref{sec:critical-euclidean-case} treat the hyperbolic and Euclidean critical cases, respectively, and present applications.

\medskip  
\textit{\textbf{Notation.}} 
Throughout the paper, $c, C>0$ denote constants whose values may change from line to line, with their dependence on parameters indicated in parentheses when needed. 
For $a,b \in \mathbb{R}$, we write $a \vee b:= \max{\{a, b\}}$ and $(a-b)_+:= (a-b) \vee 0$.
$\mathbb{N}_0=\{0,1,2,\ldots\}$

\section{Preliminaries} \label{sec:preliminaries}
In this section, we fix the notation for the underlying random walk and the branching random walk, and record some basic facts that will be used throughout the paper.

\subsection{Random walks on groups} \label{subsec:rw}
Let $\Gamma$ be a finitely generated group with identity $e$, and let $\mu$ be an admissible symmetric probability measure on $\Gamma$, where admissibility means that the support of $\mu$ generates $\Gamma$ as a semigroup.
The random walk $(Z_k)_{k\geq 0}$ on $\Gamma$ with step distribution $\mu$ is a  Markov chain with transition probability
\begin{equation*}
	\mathbf{P}(Z_{k+1}=y \mid Z_k=x) =p(x,y) =\mu(x^{-1}y), \qquad x,y\in\Gamma,\quad k\geq0.
\end{equation*}
We denote by $\mathbf{P}_x$ the law  of the RW started at $x$, and write $\mathbf{P}:=\mathbf{P}_e$. 

For $x, y\in \Gamma$, define the (weighted) Green function by 
\begin{equation} \label{eq:def-Green-function}
	G_s(x, y) :=\sum_{k=0}^{\infty} s^k p_k(x,y), \quad s\geq 0
\end{equation}
where $p_k(x,y):=\mathbf{P}_x(Z_k=y)=p_k(e,x^{-1}y)$ is the $k$-step transition probability. 
Let
\begin{equation*}
     \rho:=\limsup_{n\to\infty}p_n(x,y)^{1/n}\in(0,1] 
\end{equation*}
be the spectral radius of the RW, and set $R:=\rho^{-1} \geq 1$. 
Then $R$ is the radius of convergence of the Green function. 
Since $\mu$ is admissible, $\rho$ and $R$ are independent of $x$ and $y$; see \cite[Lemma~1.7]{Woe00}.

We have $G_s(x,y)<\infty$ for all $0<s<R$ and all $x,y\in\Gamma$, whereas finiteness at $s=R$ need not hold in general.
For $s>0$, whenever the Green function is finite, define the modified Green function by
\begin{equation} \label{eq:def-modified-Green-function} 
     \widetilde{G}_s(x,y):=\sum_{n=1}^{\infty}s^{n-1}p_n(x,y)=\frac{G_s(x,y)-\mathbf{1}_{\{x=y\}}}{s}. 
\end{equation}
Both $G_s$ and $\widetilde{G}_s$ are symmetric and left invariant.
Since $\mu$ is symmetric, $p_2(e,e)=\sum_{g\in\Gamma} \mu(g)^2>0$, which implies that $G_s(e,e)>1$. Therefore, for all $x,y\in\Gamma$,
\begin{equation}  \label{eq:modified-Green-comparison}
     c_sG_s(x,y) \leq \widetilde{G}_s(x,y) \leq s^{-1}G_s(x,y),
\end{equation}
where 
\begin{equation*}
    c_s:=\frac{G_s(e,e)-1}{sG_s(e,e)}>0.
\end{equation*}
Moreover, the upper bound in \eqref{eq:modified-Green-comparison} is an equality whenever $x\neq y$.

We shall repeatedly use the following standard heat kernel bound.
\begin{lemma} [{\cite[Proposition 6.6]{LP17}}] \label{lem:heat-kernel-estimate}
For all $x,y\in\Gamma$ and all integers $n\geq0$, $p_n(x,y)\leq\rho^n$.
\end{lemma}

\subsection{Branching random walks on groups}
Let $\Gamma$ be a finitely generated group with identity $e$, and let $\mu$ be an admissible symmetric probability measure on $\Gamma$.
Let $\nu$ be a probability distribution on $\mathbb{N}_0=\{0,1,2,\ldots\}$ satisfying  Assumption~\ref{ass:offspring} with mean $r:=\sum_{k=0}^{\infty}k\nu(k)$. 

A branching random walk with step distribution $\mu$ and offspring distribution $\nu$ can be realized as a tree-indexed random walk. 
Let $\mathbf{T}$ be a Galton--Watson tree with offspring distribution $\nu$, rooted at $\varnothing$.
For $u\in V(\mathbf{T})$, let $|u|$ denote the graph distance from $u$ to $\varnothing$, and let $u_i$, $0\leq i \leq |u|$, be the ancestor of $u$ in the $i$-th generation. 
Conditionally on $\mathbf{T}$, assign to each non-root vertex $u\in V(\mathbf{T})$ an independent $\Gamma$-valued random variable $\xi_u$ with common distribution $\mu$.
Then the position of $u$ in the BRW started at $a\in\Gamma$ is
\begin{equation*}
	 X(u):=a \xi_{u_1}\cdots\xi_{u_{|u|}}, \qquad X(\varnothing):=a.
\end{equation*}
We denote by $\pr_a$ and $\E_a$ the probability law and expectation of the BRW started at $a$, respectively, and write $\pr:=\pr_e$ and $\E:=\E_e$.

For $j\geq 0$, define the $j$-th binomial moment of the offspring distribution $\nu$ by
\begin{equation} \label{eq:offspring-binomial-moments}
    b_j:=\sum_{n=j}^{\infty}\binom{n}{j}\nu(n).
\end{equation}
Then $b_0=1$ and $b_1=r$. 
Moreover, for every $j\geq 1$, $b_j<\infty$ if and only if $\nu$ has finite $j$-th moment. 
In particular, under Assumption~\ref{ass:offspring}, $b_j<\infty$ for every $j\geq 0$.

For $x\in\Gamma$, define the local time at $x$ by
\begin{equation} \label{eq:def-local-time}
     L(x):=\sum_{u\in V(\mathbf{T})}\mathbf{1}_{\{X(u)=x\}}.	
\end{equation}
By the many-to-one formula, $\E_a[L(x)]=G_r(a,x)$. 

A BRW is called transient if $L(x)<\infty$ a.s. for all $x\in \Gamma$, and recurrent otherwise. 
It is transient if and only if $r\leq R$, where $R$ is the radius of convergence of the Green function; see \cite{BP94} for the case $r\neq R$ and \cite{GM06} for the critical case $r=R$.
Throughout the paper, we work in the transient phase.

\section{Diagrammatic expansion of moments} \label{sec:diagrammatic-expansion}
In this section, we recall the diagrammatic expansion for moments of BRW local times from \cite{GL26}, which is based on the framework of \cite{AHJ21}. 
We also derive the corresponding expansion for mixed factorial moments and establish recursive bounds for the two-point diagrams.

A \textbf{rooted plane tree} is a locally finite rooted tree $S$ with a distinguished root and a total order on the children of each vertex.
We write $V(S)$ and $E(S)$ for its vertex and edge sets, respectively, and let $c(u)$ be the number of children of $u \in V(S)$.

Fix an integer $k\geq0$. 
A \textbf{$k$-labelled rooted plane tree} is a finite rooted plane tree $S$ together with a (not necessarily injective) labelling map $ \ell:\{0,1,\ldots,k\}\to V(S)$ such that $\ell(0)$ is the root of $S$ and every leaf of $S$ is labelled (i.e. belongs to $\ell(\{0,1,\dots,k\})$).
Thus a vertex may carry several labels, and internal vertices may also be labelled. 
An isomorphism of labelled rooted plane trees is an isomorphism of rooted plane trees that preserves the labelling map.

A $k$-labelled rooted plane tree $S$ is called a \textbf{$k$-skeleton} if every unlabelled vertex has at least two children, that is, $c(v)\geq 2$ for all $v\in V^\circ(S)$, where $V^{\circ}(S)$ denotes the set of unlabelled vertices.
Let $\mathcal{S}_k$ and $\mathcal{H}_k$ be sets of isomorphism class representatives for the $k$-skeletons and the $k$-labelled rooted plane trees, respectively.

For a $k$-labelled rooted plane tree $S$ and $\mathbf{x}=(x_0,x_1,\ldots,x_k)\in\Gamma^{k+1}$, define
\begin{equation*}
     \Lambda(\mathbf{x};S) := \Bigl\{\, \mathbf{y}=(y_u)_{u\in V(S)}\in \Gamma^{V(S)}:\ y_{\ell(i)}=x_i\ \text{for each }0\le i\le k\,\Bigr\}.	
\end{equation*}
In particular, $\Lambda(\mathbf{x};S)=\varnothing$ whenever $\ell(i)=\ell(j)$ but $x_i\neq x_j$. 
If $S$ is a $k$-skeleton, we define its \textbf{$S$-diagram} by
\begin{equation} \label{eq:diagram-weight}
    \mathbf{D}(\mathbf{x};S):=\sum_{\mathbf{y}\in\Lambda(\mathbf{x};S)}\prod_{\{u,v\}\in E(S)}\widetilde{G}_r(y_u,y_v),
\end{equation}
where $\widetilde{G}_r$ is the modified Green function defined in \eqref{eq:def-modified-Green-function}, and $r$ is offspring mean.
In particular, if $S$ is the $0$-skeleton, then $\mathbf{D}(\,\cdot\,;S)\equiv1$, while for the two-vertex $1$-skeleton, $\mathbf{D}(x_0,x_1;S)=\widetilde{G}_r(x_0,x_1)$.

For integers $p,q\geq 0$, $a,x,y\in \Gamma$, and a $(p+q)$-skeleton $S$, write
\begin{equation*}
	\D_{p,q}(a;x,y;S):=  \D \big( a, \underbrace{x,\ldots,x}_{p\text{ times}}, \underbrace{y,\ldots,y}_{q\text{ times}}; S \big). 
\end{equation*}
The following expansion of mixed moment follows from  \cite[Proposition~3.1]{GL26} when $p+q\geq 1$, while the case $p=q=0$ follows directly from the definitions. 

\begin{proposition}  \label{prop:diagrammatic-moment-expansion}
Assume that $G_r(x,y)<\infty$ for all $x,y\in\Gamma$. 
Then, for all $a,x,y \in \Gamma$ and all integers $p,q \geq 0$,
\begin{equation} \label{eq:diagrammatic-moment-expansion}
    \E_a [L(x)^pL(y)^q] = \sum_{S\in\mathcal{S}_{p+q}} \D_{p,q}(a;x,y;S) \prod_{v\in V(S)}b_{c(v)}.
\end{equation}
\end{proposition}

The upper tail bound only requires control of the ordinary mixed moments, whereas the lower tail bound relies on the Paley--Zygmund inequality applied to products of binomial coefficients of local times. 
For this purpose, we derive a corresponding expansion for mixed factorial moments. 

For integers $N,n \geq 0$, let $(N)_n$ denote the falling factorial 
\begin{equation*}
    (N)_n:=N(N-1)\cdots(N-n+1),
\end{equation*}
with the conventions $(N)_0:=1$ and $(N)_n=0$ whenever $n>N$.

For $k\geq 1$, let $\mathcal{S}_k^{\star}$ be the subset of $\mathcal{S}_k$ consisting of skeletons whose labelling maps are injective on $\{1,\ldots,k\}$.
Thus, the non-root labels of a skeleton in $\mathcal{S}_k^{\star}$ correspond to distinct vertices, although one of these vertices may be the root. 
This is precisely the class arising from the distinctness constraints in factorial moments.
The following proposition gives the corresponding expansion of mixed factorial moments.
 
\begin{proposition}  \label{prop:mixed-factorial-moment-expansion}
Assume that $G_r(x,y)<\infty$ for all $x,y\in\Gamma$. 
Then, for all $a, x,y\in\Gamma$ with $x\neq y$ and all integers $p,q\geq1$,
\begin{equation}\label{eq:mixed-factorial-moment-expansion}
    \E_a\left[(L(x))_p(L(y))_q\right]  = \sum_{S\in\mathcal{S}_{p+q}^{\star}} \mathbf{D}_{p,q}(a;x,y;S) \prod_{u\in V(S)}b_{c(u)}.
\end{equation}
\end{proposition}

\begin{proof}
The proof follows the argument of \cite[Proposition~3.1]{GL26}; see also \cite[Proposition~3.1]{AHJ21}. 
The only additional point is the distinctness constraint arising from the falling factorials.  
Here we include the details for readers' convenience.

Since $x\neq y$, we have
\begin{equation*}
    (L(x))_p(L(y))_q = \#\left\{ (\mathbf{v},\mathbf{w})\in V(\mathbf{T})^p\times V(\mathbf{T})^q:  
    \begin{array}{l}
        v_1,\ldots,v_p\text{ are pairwise distinct},\\
        w_1,\ldots,w_q\text{ are pairwise distinct},\\
        X(v_i)=x,\ 1\leq i\leq p,\\
        X(w_j)=y,\ 1\leq j\leq q
    \end{array} \right\}.
\end{equation*}	
Moreover, the sets $\{ v_1,\ldots,v_p \}$ and $ \{w_1,\ldots,w_q \}$ are disjoint.
We regard $\mathbf{T}$ as a rooted plane tree by assigning an order to the children of each vertex.

For each such pair $(\mathbf{v},\mathbf{w})$, let $H(\mathbf{v},\mathbf{w})$ be the labelled rooted plane tree spanned by the root $\varnothing$ together with $v_1,\ldots,v_p,w_1,\ldots,w_q$, with labelling
\begin{equation*}
    \ell(0)=\varnothing,\qquad \ell(i)=v_i,\quad 1\leq i\leq p,\qquad
     \text{and} \qquad \ell(p+j)=w_j, \quad 1\leq j\leq q.
\end{equation*}
The distinctness of the selected vertices implies that $\ell$ is injective on $\{1,\ldots,p+q\}$. 
Let $\mathcal{H}_{p+q}^{\star}$ be the subset of $\mathcal{H}_{p+q}$ consisting of trees whose labelling maps are injective on $\{1,\ldots,p+q\}$.
Then
\begin{equation*}
    (L(x))_p(L(y))_q = \sum_{H\in\mathcal{H}_{p+q}^{\star}} \#\left\{  (\mathbf{v},\mathbf{w}): H(\mathbf{v},\mathbf{w})\cong H,\ X(v_i)=x,\ X(w_j)=y \right\}.
\end{equation*}

Set $\boldsymbol{\chi}=(a,\underbrace{x,\ldots,x}_{p\text{ times}},\underbrace{y,\ldots,y}_{q\text{ times}})$. 
Conditional on the Galton--Watson tree $\mathbf{T}$ and using the independence of the random walk increments along its edges, we have
\begin{equation*}
     \E_a\left[(L(x))_p(L(y))_q \mid \mathbf{T}\right]  =\sum_{H\in\mathcal{H}_{p+q}^{\star}} \#\{\text{embeddings of }H\text{ into }\mathbf{T}\} \sum_{\mathbf{z} \in\Lambda(\boldsymbol{\chi};H)}
    \prod_{\{u,v\}\in E(H)} p(z_u,z_v).
\end{equation*}
Here embeddings are injective maps that preserve the root, edges, and the order of children.
Taking expectation over $\mathbf{T}$ and using the standard embedding--count identity $\E\left[\#\{\text{embeddings of }H\text{ into }\mathbf{T}\}\right] = \prod_{u\in V(H)}b_{c(u)}$ (see \cite[Eq.~(3.1)]{AHJ21}), we obtain
\begin{equation} \label{eq:factorial-expansion-before-contraction}
     \E_a\left[(L(x))_p(L(y))_q\right] = \sum_{H\in\mathcal{H}_{p+q}^{\star}} \sum_{\mathbf{z}\in\Lambda(\boldsymbol{\chi};H)} \left(\prod_{u\in V(H)}b_{c(u)}\right) \left(\prod_{\{u,v\}\in E(H)}p(z_u,z_v)\right).
\end{equation}

For $H\in\mathcal{H}_{p+q}^{\star}$, let $S(H)$ be the $(p+q)$-skeleton obtained by contracting every path whose interior vertices are unlabelled vertices with exactly one child into a single edge. 
Since no labelled vertices are removed or identified in this operation, the labelling map of $S(H)$ remains injective on $\{1,\ldots,p+q\}$, and hence $S(H)\in\mathcal{S}_{p+q}^{\star}$.
On the other hand, for every $S\in\mathcal{S}_{p+q}^{\star}$, the trees $H \in\mathcal{H}_{p+q}^{\star}$ satisfying $S(H)=S$ are precisely those obtained from $S$ by replacing each edge by a path of arbitrary positive integer length.

Fix an edge $\{u,v\}\in E(S)$. 
If this edge is replaced by a path of length $n\geq1$, then summing over the spatial positions of the $n-1$ newly inserted vertices yields the $n$-step transition probability $p_n(z_u,z_v)$. 
Each newly inserted vertex is unlabelled and has exactly one child, and hence contributes the factor $b_1=r$. 
Thus, a path of length $n$ contributes $r^{n-1} p_n(z_u,z_v)$. 
Summing over all possible lengths $n\geq1$, we get
\begin{equation*}
    \sum_{n=1}^{\infty}r^{n-1}p_n(z_u, z_v) = \widetilde{G}_r(z_u, z_v).
\end{equation*}
Consequently, for every $S\in\mathcal{S}_{p+q}^{\star}$, we have
\begin{align*}
   &\sum_{\substack{H\in\mathcal{H}_{p+q}^{\star}\\S(H)=S}} \sum_{\mathbf{z}\in\Lambda(\boldsymbol{\chi};H)} \left(\prod_{u\in V(H)}b_{c(u)}\right)  \left(\prod_{\{u,v\}\in E(H)}p(z_u, z_v)\right) \\
    &\quad= \sum_{\mathbf{z}\in\Lambda(\boldsymbol{\chi};S)}  \left(\prod_{u\in V(S)}b_{c(u)}\right) \left(\prod_{\{u,v\}\in E(S)}  \widetilde{G}_r(z_u, z_v)\right) = \D_{p,q}(a;x,y;S)  \prod_{u\in V(S)}b_{c(u)}.
\end{align*}
Therefore, the desired result follows by summing over $S\in\mathcal{S}_{p+q}^{\star}$ in \eqref{eq:factorial-expansion-before-contraction}.
\end{proof}

The expansions in Propositions~\ref{prop:diagrammatic-moment-expansion} and \ref{prop:mixed-factorial-moment-expansion} reduce moment estimates to bounds on diagrams. 
Instead of estimating each skeleton separately, we introduce the notion of maximal diagrams. 

For $k\geq 0$, let $\mathcal{S}_k'$ be the subset of $\mathcal{S}_k$ consisting of skeletons whose labelling maps are injective on $\{0,1,\ldots,k\}$.
Then, $\mathcal{S}_k' \subseteq \mathcal{S}_k^{\star} \subseteq \mathcal{S}_k$ for all $k\geq1$. 
For $a,x,y\in\Gamma$ and integers $p,q \geq 0$ with $p+q\geq 1$, define
\begin{equation}\label{eq:def-maximal-diagram}
    M_{p,q}(a;x,y):=\max_{S\in\mathcal{S}_{p+q}'}\D_{p,q}(a;x,y;S),
\end{equation}
and set $M_{0,0}(a;x,y):=\mathbf{1}_{ \{a=x=y\} }$. 
Note that 
\begin{equation} \label{eq:one-point-max}
    M_{1,0}(a;x,y)=\widetilde{G}_r(a,x) \quad\text{and}\quad  M_{0,1}(a;x,y)=\widetilde{G}_r(a,y).
\end{equation}

The following proposition provides a recursive bound for maximal diagrams.

\begin{proposition}  \label{prop:maximal-diagram-recursion}
Assume that $G_r(x,y)<\infty$ for all $x,y\in\Gamma$. 
Then for all $a,x,y\in\Gamma$ and all integers $p,q\geq0$ with $p+q\geq 2$,
\begin{equation} \label{eq:maximal-diagram-recursion}
    M_{p,q}(a;x,y) \leq C_0 \max_{\substack{0\leq p_1\leq p,\ 0\leq q_1\leq q\\1\leq p_1+q_1\leq p+q-1}} \sum_{z\in\Gamma} \widetilde{G}_r(a,z)  M_{p_1,q_1}(z;x,y) M_{p-p_1,q-q_1}(z;x,y),
\end{equation}
where $C_0:=1 \vee \widetilde{G}_r(e,e)^{-1}$.
\end{proposition}

\begin{proof}
Fix $S\in\mathcal{S}_{p+q}'$, and let $\ell:\{0,1,\ldots,p+q\}\to V(S)$ be its labelling map. 
Since $\ell$ is injective, the root carries only the root label $0$. 
For $u\in V(S)$, let $S_u$ denote the descendant subtree rooted at $u$, including $u$ itself.
We distinguish two cases according to the number of children of the root.

\medskip
\noindent\textbf{Case 1. The root has at least two children.}
Choose one child $u$ of the root and set
\begin{equation*}
    S_1:=S_u\cup\{\varnothing\} \qquad \text{and} \qquad S_2:=S\setminus S_u.
\end{equation*}
Since every leaf of $S$ is labelled and the root carries only label $0$, both $S_1$ and $S_2$ contain at least one non-root label.
We regard $S_1$ and $S_2$ as labelled skeletons rooted at $\varnothing$, assigning the root label $0$ and relabelling the inherited non-root labels in their original order.

Suppose that $S_1$ contains $p_1$ labels assigned to $x$ and $q_1$ labels assigned to $y$. 
Then $1\leq p_1+q_1\leq p+q-1$, and $S_2$ contains the remaining $p-p_1$ labels assigned to $x$ and $q-q_1$ labels assigned to $y$.
By factorization of the diagram at the root, we obtain
\begin{align*}
     \D_{p,q}(a;x,y;S) &= \D_{p_1,q_1}(a;x,y;S_1) \D_{p-p_1,q-q_1}(a;x,y;S_2) \\
	 &\leq M_{p_1,q_1}(a;x,y) M_{p-p_1,q-q_1}(a;x,y).
\end{align*}

\medskip
\noindent\textbf{Case 2. The root has exactly one child.}
Let $w$ be this child. 
If $w$ is assigned the spatial position $z$, then the edge from the root to $w$ contributes the factor $\widetilde{G}_r(a,z)$. 
We distinguish according to whether $w$ carries a non-root label.

\smallskip
\noindent\textbf{Case 2a. The vertex $w$ carries no non-root label.}
Since $S$ is a skeleton, the vertex $w$ has at least two children. 
Choose one child $u$ of $w$ and set
\begin{equation*}
    S_1:=S_u\cup\{w\} \qquad \text{and} \qquad S_2:=S_w\setminus S_u.
\end{equation*}
We regard $S_1$ and $S_2$ as labelled skeletons rooted at $w$, assigning the new root label $0$ to $w$ and relabelling the inherited non-root labels in their original order.

Suppose that $S_1$ contains $p_1$ labels assigned to $x$ and $q_1$ labels assigned to $y$. 
Then $1\leq p_1+q_1\leq p+q-1$.
Factoring the diagram at $w$ yields
\begin{align*}
      \D_{p,q}(a;x,y;S)	 
    &=\sum_{z\in \Gamma}\widetilde{G}_r(a,z) \ \D_{p_1,q_1}(z;x,y;S_1) \ \D_{p-p_1,q-q_1}(z;x,y;S_2) \\ 
    &\leq \sum_{z\in \Gamma} \widetilde{G}_r(a,z) \ M_{p_1,q_1}(z;x,y) \ M_{p-p_1,q-q_1}(z;x,y).
\end{align*}

\smallskip
\noindent\textbf{Case 2b. The vertex $w$ carries a non-root label.}
By injectivity of $\ell$, the vertex $w$ carries exactly one such label. 
Suppose first that this label is assigned to $x$. 
Then the spatial position of $w$ is forced to be $x$. 
We restrict to the subtree $S_w$, remove the label at $w$, assign $w$ the new root label $0$, and relabel the non-root labels on the other vertices of $S_w$ in their original order.
The resulting labelled skeleton, denoted by $S_w^{\circ}$, contains $p-1$ labels assigned to $x$ and $q$ labels assigned to $y$. 
Hence
\begin{equation*}
    \D_{p,q}(a;x,y;S) = \widetilde{G}_r(a,x)\D_{p-1,q}(x;x,y;S_w^{\circ}) \leq \widetilde{G}_r(a,x)M_{p-1,q}(x;x,y).
\end{equation*}
Since $M_{1,0}(x;x,y)=\widetilde{G}_r(x,x)=\widetilde{G}_r(e,e)$, we obtain
\begin{equation*}
     \D_{p,q}(a;x,y;S) \leq \widetilde{G}_r(e,e)^{-1} \sum_{z\in\Gamma}\widetilde{G}_r(a,z)  M_{1,0}(z;x,y)M_{p-1,q}(z;x,y).
\end{equation*}
The case where $w$ carries a label assigned to $y$ follows analogously, using $M_{0,1}(y;x,y)=\widetilde{G}_r(e,e)$.

Combining the above cases and taking the maximum over $S\in\mathcal{S}_{p+q}'$,  it remains to absorb the contribution from Case 1.
Since $\widetilde{G}_r(a,a)=\widetilde{G}_r(e,e)$, this contribution is bounded by the convolution term in \eqref{eq:maximal-diagram-recursion} by retaining only the term $z=a$, at the cost of the factor $\widetilde{G}_r(e,e)^{-1} \leq C_0$. 
This completes the proof of \eqref{eq:maximal-diagram-recursion}.
\end{proof}

We now extend this bound to all labelling map.
By reducing a general diagram to one with an injective labelling map, we can control all diagrams.

\begin{lemma} \label{lem:noninjective-diagram-reduction}
For all $a,x,y\in\Gamma$ and all integers $p,q\geq1$,
\begin{align*}
	&\max_{S\in \mathcal{S}_{p+q}} \D_{p,q}(a;x,y;S)  \\
	&\quad \leq  \max\Bigg\{
        \max_{\substack{1\le i\le p\\1\le j\le q}}M_{i,j}(a;x,y),
        \max_{1\le i\le p}M_{i,0}(a;x,y),
        \max_{1\le j\le q}M_{0,j}(a;x,y),
        M_{0,0}(a;x,y)
    \Bigg\}.
\end{align*}
Moreover, the term involving $M_{p',0}$ can occur only if $y=a$ or $y=x$, the term involving $M_{0,q'}$ only if $x=a$, and the term $M_{0,0}$ only if $y=x=a$.
\end{lemma}

\begin{proof}
Fix $S\in\mathcal{S}_{p+q}$ and suppose that $\D_{p,q}(a;x,y;S)>0$. 
At each labelled vertex, retain only the label with the smallest index and delete all the others.
Since the root label is $0$, it is always retained.
Let $p'$ and $q'$ be the numbers of retained labels assigned to $x$ and $y$,  respectively, and let $S'$ be the resulting labelled rooted plane tree after relabelling the retained non-root labels in their original order.

Since the diagram is positive, labels attached to the same vertex must be assigned the same spatial position.  
Hence deleting repeated labels does not change the admissible spatial assignments or the edge set.
Moreover, every previously labelled vertex remains labelled, so $S'$ is again a skeleton.
Therefore, 
\begin{equation}\label{eq:diagram-label-reduction}
    \D_{p,q}(a;x,y;S) = \D_{p',q'}(a;x,y;S').
\end{equation}
By construction, the labelling of $S'$ is injective. Thus,
\begin{align*}
	&\D_{p,q}(a;x,y;S)  \\
	&\quad \leq  \max\Bigg\{
        \max_{\substack{1\le i\le p\\1\le j\le q}}M_{i,j}(a;x,y),
        \max_{1\le i\le p}M_{i,0}(a;x,y),
        \max_{1\le j\le q}M_{0,j}(a;x,y),
        M_{0,0}(a;x,y)
    \Bigg\}.
\end{align*}
Taking the maximum over $S\in\mathcal{S}_{p+q}$, we get the desired inequality.

We now identify the conditions under which the reduction can yield the last three types of diagrams.
Because the labels assigned to $x$ have smaller indices than those assigned to $y$, $q'=0$ can occur only if every label assigned to $y$ shares its vertex with either the root label or a label assigned to $x$, which requires $y=a$ or $y=x$. 
Similarly, $p'=0$ requires $x=a$, while $p'=q'=0$ requires $a=x=y$. 
This completes the proof.
\end{proof}

\section{A general criterion for two-point tails} \label{sec:general-two-point-criterion}
In this section, we first prove Theorem~\ref{thm:two-point-tail-estimate}.
We then verify the corresponding conditions on the Green function in the subcritical case in order to obtain  Corollary~\ref{cor:subcritical-two-point-tail-estimate}.

\subsection{Consequences of the Green function criterion}
We first record some consequences of the  conditions \eqref{eq:criterion-triangle} and \eqref{eq:criterion-mixed} imposed on the Green function.
By symmetry, the triangle condition \eqref{eq:criterion-triangle} also holds with the square on $G_r(z,x)$.
Hence, by the Cauchy--Schwarz inequality, we obtian
\begin{align}
    \Phi_r(a;x,y)^2  
    &\leq \left(\sum_{z\in\Gamma}G_r(a,z)G_r(z,x)^2\right)\left(\sum_{z\in\Gamma}G_r(a,z)G_r(z,y)^2\right) \notag\\
    &\leq C G_r(a,x)G_r(a,y). \label{eq:criterion-Phi-CS}
\end{align}
In particular, $\Phi_r(a;x,y)<\infty$. 
Applying this bound with $a=z$ and summing against $G_r(a,z)$, we have the following estimate.

\begin{proposition} \label{prop:Phi-square-convolution}
Under the assumptions of Theorem~\ref{thm:two-point-tail-estimate}, there exists a constant $C>0$ such that, for all $a,x,y\in\Gamma$,
\begin{equation*} 
    \sum_{z\in\Gamma} G_r(a,z)\Phi_r(z;x,y)^2  \leq C \Phi_r(a;x,y).
\end{equation*}
\end{proposition}

Combining Proposition~\ref{prop:Phi-square-convolution}, the conditions  \eqref{eq:criterion-triangle} and \eqref{eq:criterion-mixed}, the definition of $\Phi_r$, and  $\widetilde{G}_r\leq r^{-1}G_r$, we immediately have the following convolution bounds needed in the two-point recursion.

\begin{corollary} \label{cor:modified-Green-function-estimates}
Under the assumptions of Theorem~\ref{thm:two-point-tail-estimate}, there exists a constant $C>0$ such that, for all $a,x,y\in\Gamma$,
\begin{align}
    \sum_{z\in\Gamma}\widetilde{G}_r(a,z)G_r(z,x)G_r(z,y)
    &\leq C\Phi_r(a;x,y), \label{eq:convolution-1}\\
    \sum_{z\in\Gamma}\widetilde{G}_r(a,z)G_r(z,x)\Phi_r(z;x,y)
    &\leq C\Phi_r(a;x,y), \label{eq:convolution-2}\\
    \sum_{z\in\Gamma}\widetilde{G}_r(a,z)\Phi_r(z;x,y)^2
    &\leq C\Phi_r(a;x,y). \label{eq:convolution-3}
\end{align}
The second estimate also holds with $x$ and $y$ interchanged.
\end{corollary}

\subsection{Moment estimates under the criterion} \label{subsec:general-two-point-moments}

We now apply Proposition~\ref{prop:maximal-diagram-recursion} with Corollary~\ref{cor:modified-Green-function-estimates} to obtain moment estimates.
We begin with the one-point maximal diagram estimate, corresponding to the case where  one of the two groups of labels is empty. 
It is obtained by specializing Proposition~\ref{prop:maximal-diagram-recursion} to a single target, and the resulting induction is analogous to that used to prove  \cite[(4.2)]{GL26}.

\begin{lemma}  \label{lem:one-point-diagram-bound}
Under the assumptions of Theorem~\ref{thm:two-point-tail-estimate}, there exists $\eta_0\geq 1$ such that, for all $a,x,y\in\Gamma$ and all integers $p,q\geq1$,
\begin{equation} \label{eq:one-point-maximal-diagram-bound}
    M_{p,0}(a;x,y)\leq r^{-1}\eta_0^{p-1}G_r(a,x)\qquad \text{and} \qquad
    M_{0,q}(a;x,y)\leq r^{-1}\eta_0^{q-1}G_r(a,y).
\end{equation}
\end{lemma}

\begin{proof}
We prove the first bound and the second follows similarly. For $p=1$, by \eqref{eq:one-point-max}, we have
\begin{equation*}
    M_{1,0}(a;x,y)=\widetilde{G}_r(a,x) \leq r^{-1}G_r(a,x).
\end{equation*}
Suppose that the bound holds for all positive integers less than $p$. Then by the one-point specialization of Proposition~\ref{prop:maximal-diagram-recursion}, we get 
\begin{align*}
    M_{p,0}(a;x,y)
    &\leq C_0\max_{1\leq i\leq p-1} \sum_{z\in\Gamma}\widetilde G_r(a,z) M_{i,0}(z;x,y)M_{p-i,0}(z;x,y) \\
    &\leq C_0r^{-3}\eta_0^{p-2} \sum_{z\in\Gamma}G_r(a,z)G_r(z,x)^2 \\
    &\leq C_0Cr^{-3}\eta_0^{p-2}G_r(a,x),
\end{align*}
where the last inequality follows from the triangle condition~\eqref{eq:criterion-triangle}. 
Choosing $\eta_0\geq C_0Cr^{-2} \vee 1$, we complete the induction.
\end{proof}

We next bound the two-point maximal diagrams.
\begin{proposition} \label{prop:two-point-maximal-diagram-bound}
Under the assumptions of Theorem~\ref{thm:two-point-tail-estimate}, there exist constants $C<\infty$ and $\eta\geq1$ such that, for all $a,x,y\in\Gamma$ and all integers $p,q\geq1$,
\begin{equation} \label{eq:two-point-maximal-diagram-bound}
    M_{p,q}(a;x,y)\leq C\eta^{p+q} \Phi_r(a;x,y).
\end{equation}
\end{proposition}

\begin{proof}
We prove the stronger estimate
\begin{equation} \label{eq:strong-two-point-maximal-diagram-bound}
    M_{p,q}(a;x,y)\leq A \eta^{p+q-2} \Phi_r(a;x,y),
\end{equation}
for suitable constants $A<\infty$ and $\eta\geq\eta_0$.

Let $C_1<\infty$ be a constant for which all the estimates in Corollary~\ref{cor:modified-Green-function-estimates} hold. 
We argue by induction on $N=p+q$. 
For $N=2$, we necessarily have $p=q=1$. Then by Proposition~\ref{prop:maximal-diagram-recursion} and Lemma~\ref{lem:one-point-diagram-bound},
\begin{align*}
    M_{1,1}(a;x,y)
    &\leq C_0\sum_{z\in\Gamma}\widetilde{G}_r(a,z) M_{1,0}(z;x,y)M_{0,1}(z;x,y)\\
    &\leq C_0 r^{-2}\sum_{z\in\Gamma}\widetilde{G}_r(a,z)  G_r(z,x)G_r(z,y)\\
    &\leq C_0C_1 r^{-2}\Phi_r(a;x,y).
\end{align*}
Hence, \eqref{eq:strong-two-point-maximal-diagram-bound} holds for $N=2$ provided that $A\geq C_0C_1r^{-2}$.

Suppose now that $N\geq 3$ and that \eqref{eq:strong-two-point-maximal-diagram-bound} holds for all positive integers $i,j \geq 1$ with $i+j<N$. 
Let $p,q\geq1$ satisfy $p+q=N$. 
For a pair $(p_1,q_1)$ appearing in Proposition~\ref{prop:maximal-diagram-recursion}, set $ p_2:=p-p_1$ and $q_2:=q-q_1$. 
Both $(p_1,q_1)$ and $(p_2,q_2)$ are nonzero. 
According to whether these two pairs have one or two positive coordinates, there are three cases.

\smallskip
\noindent\textbf{Case 1. Both factors are one-point diagrams.}
Since $p,q\geq 1$, one factor contains only labels assigned to $x$ and the other only labels assigned to $y$. 
Hence, up to interchanging the two factors, we may assume that 
\begin{equation*}
    (p_1,q_1)=(p,0) \qquad \text{and} \qquad (p_2,q_2)=(0,q).
\end{equation*}
Consequently, by Lemma~\ref{lem:one-point-diagram-bound} and Corollary~\ref{cor:modified-Green-function-estimates}, we have
\begin{align*}
	 \sum_{z\in\Gamma} \widetilde{G}_r(a,z) M_{p,0}(z;x,y)M_{0,q}(z;x,y)   
	  &\leq r^{-2} \eta^{p+q-2} \sum_{z\in\Gamma} \widetilde{G}_r(a,z)G_r(z,x)G_r(z,y) \\
	  &\leq C_1 r^{-2} \eta^{p+q-2}\Phi_r(a;x,y).
\end{align*}

\smallskip
\noindent\textbf{Case 2. One factor is a one-point diagram and the other is a two-point diagram.}
By the symmetry between the two factors, we may assume that the first factor is a one-point diagram.
We first consider the case in which it contains only labels assigned to $x$, while the second factor contains labels assigned to both $x$ and $y$. 
Thus,
\begin{equation*}
    (p_1,q_1)=(i,0) \qquad \text{and} \qquad  (p_2,q_2)=(p-i,q)
\end{equation*}
for some $1\leq i\leq p-1$. 
Then, by Lemma~\ref{lem:one-point-diagram-bound}, the induction hypothesis, and Corollary~\ref{cor:modified-Green-function-estimates}, we obtain
\begin{align*}
	 \sum_{z\in\Gamma} \widetilde{G}_r(a,z)  M_{i,0}(z;x,y) M_{p-i,q}(z;x,y) &\leq A r^{-1}\eta^{p+q-3}  \sum_{z\in\Gamma} \widetilde{G}_r(a,z)G_r(z,x)\Phi_r(z;x,y)\\
	 &\leq C_1 A r^{-1} \eta^{p+q-3} \Phi_r(a;x,y).
\end{align*}
The case in which the one-point factor contains only labels assigned to $y$ follows analogously by interchanging $x$ and $y$ in the second estimate of Corollary~\ref{cor:modified-Green-function-estimates}.

\smallskip
\noindent\textbf{Case 3. Both factors are two-point diagrams.}
In this case, both factors contain labels assigned to both $x$ and $y$, so $p_1,q_1,p_2,q_2 \geq 1$.  
Thus, by the induction hypothesis and Corollary~\ref{cor:modified-Green-function-estimates}, we obtain
\begin{align*}
	  \sum_{z\in\Gamma} \widetilde{G}_r(a,z)  M_{p_1,q_1}(z;x,y) M_{p_2,q_2}(z;x,y)
	  &\leq A^2 \eta^{p+q-4}  \sum_{z\in\Gamma} \widetilde{G}_r(a,z) \Phi_r(z;x,y)^2 \\
	  &\leq C_1 A^2 \eta^{p+q-4} \Phi_r(a;x,y).
\end{align*}

Since the three cases above cover all possible pairs $(p_1, q_1)$, by Proposition~\ref{prop:maximal-diagram-recursion}, we have
\begin{equation*}
\begin{aligned}
    M_{p,q}(a;x,y)
    &\leq C_0C_1r^{-2}\eta^{p+q-2}\Phi_r(a;x,y) \ \vee\ C_0C_1Ar^{-1}\eta^{p+q-3}\Phi_r(a;x,y)\\
    &\quad\vee\ C_0C_1A^2\eta^{p+q-4}\Phi_r(a;x,y).
\end{aligned}
\end{equation*}
Choose $A\geq \max\{C_0C_1r^{-2},  1\} $, and then choose $\eta \geq \eta_0$ sufficiently large so that $\eta\geq C_0C_1r^{-1}$ and $ \eta^2\geq C_0C_1A$. 
Consequently, each term is bounded by $A \eta^{p+q-2} \Phi_r(a;x,y)$. 
This completes the induction. 
Since $\eta \geq 1 $, the stated estimate follows immediately.
\end{proof}

By combining Proposition~\ref{prop:two-point-maximal-diagram-bound} with Lemma~\ref{lem:noninjective-diagram-reduction}, we can bound the diagrams uniformly for arbitrary skeletons. 

\begin{corollary} \label{cor:uniform-diagram-bound}
Under the assumptions of Theorem~\ref{thm:two-point-tail-estimate}, there exist constants $C<\infty$ and $\eta\geq 1$ such that, for all $a,x,y\in\Gamma$ and all integers $p,q\geq1$,
\begin{equation} \label{eq:uniform-diagram-bound}
    \max_{S\in\mathcal{S}_{p+q}}  \D_{p,q}(a;x,y;S)  \leq C\eta^{p+q}\Phi_r(a;x,y).
\end{equation}
\end{corollary}

\begin{proof}
By Lemma~\ref{lem:noninjective-diagram-reduction}, it suffices to control the following four quantities:
\begin{equation*}
	 \max_{\substack{1\leq i\leq p \\ 1\leq j\leq q}}  M_{i,j}(a;x,y), \quad \max_{1\leq i\leq p}  M_{i,0}(a;x,y), \quad \max_{1\le j\le q}  M_{0,j}(a;x,y), \quad M_{0,0}(a;x,y).
\end{equation*} 

The first quantity is bounded by Proposition~\ref{prop:two-point-maximal-diagram-bound}. 
For the second quantity, Lemma~\ref{lem:noninjective-diagram-reduction} shows that $M_{i,0}$ can occur only if $y=a$ or $y=x$. In either case, $\Phi_r(a;x,y)\geq G_r(e,e)^2 G_r(a,x)$. 
Hence, by Lemma~\ref{lem:one-point-diagram-bound}, after enlarging  $\eta$ if necessary so that $\eta\geq\eta_0$, we have
\begin{equation*}
    M_{i,0}(a;x,y)\leq C \eta^{i} \Phi_r(a;x,y) \leq C \eta^{p+q} \Phi_r(a;x,y).
\end{equation*}
For the third quantity, $M_{0,j}$ can occur only if $x=a$, in which case we have $\Phi_r(a;a,y)\geq G_r(e,e)^2G_r(a,y)$. 
The desired bound follows by the same argument.
Finally, since $\Phi_r(a;a,a)\geq G_r(a,a)^3 \geq 1$, we have 
\begin{equation*}
	M_{0,0}(a;x,y)=\mathbf{1}_{\{a=x=y \}} \leq \Phi_r(a;x,y).
\end{equation*}
This completes the proof.
\end{proof}

To sum the diagrammatic expansion over all skeletons, we use an estimate for the skeleton partition function.
Since $\nu$ has a finite exponential moment, there exist constants $c,C>0$ such that $\nu(k) \leq C e^{-ck}$ for all integers $k\geq 0$. 
Therefore, the proof of \cite[Lemma~5.2]{AHJ21} applies and yields the following result.

\begin{lemma} [{\cite[Lemma~5.2]{AHJ21}}] \label{lem:the-skeleton-partition-function}
Under Assumption~\ref{ass:offspring}, there exist constants $C=C(\nu)<\infty$ and $\lambda=\lambda(\nu)\geq1$ such that, for all integers $k\geq1$,
\begin{equation} \label{eq:skeleton-partition-estimate}
    \sum_{S\in\mathcal{S}_k}\prod_{u\in V(S)}b_{c(u)} \leq C\lambda^k k!.
\end{equation}
\end{lemma}

Combining the skeleton partition estimate with the uniform diagram bound, we obtain the following mixed moment upper bound.

\begin{corollary} \label{cor:mixed-moment-upper-bound}
Under the assumptions of Theorem~\ref{thm:two-point-tail-estimate}, there exist constants $C<\infty$ and $\zeta \geq 1$ such that, for all $a,x,y\in\Gamma$ and all integers $p,q\geq1$,
\begin{equation} \label{eq:mixed-moment-upper-bound}
    \E_a\bigl[L(x)^pL(y)^q\bigr] \leq C\zeta^{p+q}p!q!\Phi_r(a;x,y).
\end{equation}
\end{corollary}

\begin{proof}
By Proposition~\ref{prop:diagrammatic-moment-expansion}, Corollary~\ref{cor:uniform-diagram-bound}, and Lemma~\ref{lem:the-skeleton-partition-function}, we obtain
\begin{align*}
    \E_a\left[L(x)^p L(y)^q\right]
    &=\sum_{S\in\mathcal{S}_{p+q}}\D_{p,q}(a;x,y;S)\prod_{u\in V(S)}b_{c(u)} \leq C\eta^{p+q} \lambda^{p+q}(p+q)! \Phi_r(a;x,y)  \\
    &\leq C(2\lambda\eta)^{p+q}p!q!  \Phi_r(a;x,y),
\end{align*}
where the last inequality follows from  $(p+q)!=\binom{p+q}{p}p!q!\leq2^{p+q}p!q!$.
\end{proof}

In order to establish the lower bound in Theorem~\ref{thm:two-point-tail-estimate}, we derive the following estimate for mixed factorial moments.

\begin{corollary} \label{cor:mixed-factorial-moment-lower-bound}
Under the assumptions of Theorem~\ref{thm:two-point-tail-estimate}, there exist constants $c,\alpha>0$ such that, for all $a,x,y\in\Gamma$ with $x\neq y$ and all integers $p,q\geq1$,
\begin{equation} \label{eq:mixed-factorial-moment-lower-bound}
    \E_a\bigl[(L(x))_p (L(y))_q\bigr] \geq c \alpha^{p+q}p!q!\Phi_r(a;x,y).
\end{equation}
\end{corollary}

\begin{proof}
In the expansion of Proposition~\ref{prop:mixed-factorial-moment-expansion}, we consider only the following family of skeletons.  
The root has a single unlabelled child $v$, and $v$ has two children whose planar order is fixed. 
The left subtree is a linear chain consisting of exactly $p$ labelled vertices with labels $1,\ldots,p$ assigned to $x$,  while the right subtree is a linear chain consisting of exactly $q$ labelled vertices with labels  $p+1,\ldots,p+q$ assigned to $y$.
The labels may be arranged arbitrarily along their respective chains, giving $p!q!$ such skeletons.

For each such skeleton,
\begin{equation*}
    \prod_{u\in V(S)}b_{c(u)}=b_2 r^{p+q-1},
\end{equation*}
and the corresponding diagram is
\begin{equation*}
    \D_{p,q}(a;x,y;S) =\widetilde G_r(e,e)^{p+q-2} \sum_{z\in\Gamma} \widetilde{G}_r(a,z)\widetilde G_r(z,x)\widetilde G_r(z,y).
\end{equation*}
Hence
\begin{align*}
    \E_a\bigl[(L(x))_p(L(y))_q\bigr]
    &\geq p!q!b_2r^{p+q-1}\widetilde G_r(e,e)^{p+q-2}  \sum_{z\in\Gamma}  \widetilde{G}_r(a,z)\widetilde{G}_r(z,x)\widetilde{G}_r(z,y) \\
    &\geq c_r^3 p!q!b_2r^{p+q-1}\widetilde{G}_r(e,e)^{p+q-2}\Phi_r(a;x,y),
\end{align*}
where the last inequality follows from \eqref{eq:modified-Green-comparison}.
Since $\nu(0)+\nu(1)<1$, we have $b_2>0$, and \eqref{eq:modified-Green-comparison} also implies $\widetilde{G}_r(e,e)>0$. 
Hence, setting $\alpha= r \widetilde{G}_r(e,e)$ and choosing $c>0$ appropriately, we obtain the desired result.
\end{proof}

\subsection{Proof of the general criterion} \label{subsec:general-two-point-tail}
In this subsection, we prove Theorem~\ref{thm:two-point-tail-estimate}.
We first establish the two-sided tail estimate and then show that the exponential constants necessarily satisfy $c_2>C_2$. 

\begin{proposition} \label{prop:two-sided-estimates}
Under the assumptions of Theorem~\ref{thm:two-point-tail-estimate}, there exist positive constants $c_1,c_2,C_1,C_2$ such that the two-sided tail estimate \eqref{eq:two-point-tail-estimate} holds for all $a,x,y \in \Gamma$ with $x\neq y$ and all integers $n,m \geq 1$ .
\end{proposition}

\begin{proof}
Let $C,\zeta,c,\alpha$ be the constants in Corollaries~\ref{cor:mixed-moment-upper-bound} and \ref{cor:mixed-factorial-moment-lower-bound}. 
By enlarging $\zeta$ if necessary, we may assume that
\begin{equation} \label{eq:zeta-alpha-comparison}
   \log \left(\frac{4\zeta^2}{\alpha^2} \right)>\frac{1}{2\zeta}.
\end{equation}

\smallskip
\noindent\textbf{Upper bound.}
By Tonelli's theorem and Corollary~\ref{cor:mixed-moment-upper-bound}, we have
\begin{align*}
    \E_a\left[\left(e^{\frac{1}{2\zeta}L(x)}-1\right)\left(e^{\frac{1}{2\zeta}L(y)}-1\right)\right]
    &=\sum_{p=1}^{\infty}\sum_{q=1}^{\infty}\frac{1}{(2\zeta)^{p+q}p!q!} \E_a[L(x)^pL(y)^q]\\
    &\leq C\Phi_r(a;x,y)\sum_{p=1}^{\infty}\sum_{q=1}^{\infty}2^{-(p+q)}\\
    &=C\Phi_r(a;x,y).
\end{align*}
On the event $\{L(x)\geq n,L(y)\geq m\}$,
\begin{equation*}
    \left(e^{\frac{1}{2\zeta}L(x)}-1\right)\left(e^{\frac{1}{2\zeta}L(y)}-1\right) \geq\left(e^{\frac{1}{2\zeta} n}-1\right)\left(e^{\frac{1}{2\zeta}m}-1\right).
\end{equation*}
Hence, by Markov's inequality, we obtain
\begin{equation*}
    \pr_a\left(L(x)\geq n,L(y)\geq m\right) \leq\frac{C\Phi_r(a;x,y)}{\left(e^{n/(2\zeta)}-1\right)\left(e^{m/(2\zeta)}-1\right)}.
\end{equation*}
Since $n\geq1$,
\begin{equation*}
    e^{\frac{n}{2\zeta}}-1 =e^{\frac{n}{2\zeta}}\left(1-e^{-\frac{n}{2\zeta} }\right)\geq \left(1-e^{-\frac{1}{2\zeta}}\right)e^{\frac{n}{2\zeta}},
\end{equation*}
and the same estimate holds with $n$ replaced by $m$. Therefore,
\begin{equation*}
     \pr_a\left(L(x)\geq n,L(y)\geq m\right) \leq \frac{C}{\left(1-e^{-1/(2\zeta)}\right)^2} e^{-\frac{1}{2\zeta}(n+m)}\Phi_r(a;x,y).
\end{equation*}

\smallskip
\noindent\textbf{Lower bound.}
For $n,m\geq1$, set
\begin{equation*}
    Z_{n,m}:=\binom{L(x)}{n}\binom{L(y)}{m} =\frac{(L(x))_n(L(y))_m}{n!m!}.
\end{equation*}
Then
\begin{equation*}
    \{Z_{n,m}>0\}=\{L(x)\geq n,L(y)\geq m\}.
\end{equation*}
Hence, by Corollary~\ref{cor:mixed-factorial-moment-lower-bound}, we have
\begin{equation}\label{eq:two-point-tail-Z-first-moment}
    \E_a[Z_{n,m}] \geq c\alpha^{n+m} \Phi_r(a;x,y).
\end{equation}
On the other hand, since $Z_{n,m}^2\leq L(x)^{2n}L(y)^{2m}/(n! m!)^2$, Corollary~\ref{cor:mixed-moment-upper-bound} gives
\begin{equation} \label{eq:two-point-tail-Z-second-moment}
    \E_a[Z_{n,m}^2] \leq C\zeta^{2(n+m)}\frac{(2n)!(2m)!}{(n!)^2(m!)^2}\Phi_r(a;x,y) \leq C(4\zeta^2)^{n+m}\Phi_r(a;x,y), 
\end{equation}
where we used  $(2k)!/(k!)^2\leq4^k$.  
Hence, by the Paley--Zygmund inequality and \eqref{eq:two-point-tail-Z-first-moment}--\eqref{eq:two-point-tail-Z-second-moment}, we obtain
\begin{equation*}
    \pr_a\left(L(x)\geq n,L(y)\geq m\right) =\pr_a(Z_{n,m}>0) \geq\frac{\E_a[Z_{n,m}]^2}{\E_a[Z_{n,m}^2]} \geq \frac{c^2}{C}\left(\frac{4\zeta^2}{\alpha^2}\right)^{-(n+m)}\Phi_r(a;x,y).
\end{equation*}
Together with \eqref{eq:zeta-alpha-comparison}, this proves  the desired two-sided estimate with $c_2= \log\frac{4\zeta^2}{\alpha^2}>\frac{1}{2\zeta}=C_2 $.
\end{proof}

The proof of Proposition~\ref{prop:two-sided-estimates} gives exponential constants $c_2$ and $C_2$ satisfying $c_2>C_2$.
We now show that this strict inequality is necessary for any choice of exponential constants in the two-sided estimate \eqref{eq:two-point-tail-estimate}.

\begin{proposition} \label{prop:strict-exponential-constants}
Under the assumptions of Theorem~\ref{thm:two-point-tail-estimate}, any constants $c_1,c_2,C_1,C_2>0$ for which \eqref{eq:two-point-tail-estimate} holds necessarily satisfy $c_2>C_2$.
\end{proposition}

\begin{proof}
We first show that $c_2 \geq C_2$. 
Fix $x\neq y$. Since $0<\Phi_r(e;x,y)<\infty$, the two-sided estimate \eqref{eq:two-point-tail-estimate} gives $c_1/ C_1 \leq e^{(c_2-C_2)(n+m)}$ for all integers $n,m\geq 1$.
Hence, letting $n+m \to\infty$, we obtain $c_2 \geq C_2$.

Now suppose, for contradiction, that $c_2=C_2=:c^*$.	
Then by \eqref{eq:two-point-tail-estimate}, there exist constants $c, C>0$ such that, for all integers $n, m\geq1$,
\begin{equation}\label{eq:equal-exponential-rates}
    c\Phi_r(e;x,y) e^{-c^*(n+m)}  \leq  \mathbb{P} (L(x)\geq n,L(y)\geq m) \leq  C\Phi_r(e;x,y) e^{-c^*(n+m)}.
\end{equation}
For $s>1$, define
\begin{equation*}
	F(s):= \mathbb{E} \left[ \left(s^{L(x)}-1 \right)   \left(s^{L(y)}-1 \right) \right].
\end{equation*}
Since $s^N-1=(s-1)\sum_{k=1}^{\infty} s^{k-1}\mathbf 1_{\{N\geq k\}}$ for all integers $N \geq 0$, by Tonelli's theorem, we have
\begin{equation}\label{eq:joint-tail-generating-function}
    F(s) = (s-1)^2 \sum_{n=1}^{\infty} \sum_{m=1}^{\infty} s^{n+m-2}  \mathbb{P}(L(x)\geq n,L(y)\geq m).
\end{equation}
Hence, the upper bound in \eqref{eq:equal-exponential-rates} implies that $F(s)<\infty$ for all $ 1<s<e^{c^*}$.

Since $\nu(0)+\nu(1) <1$, there exists an integer $k\geq2$ such that $\nu(k)>0$.
Furthermore, choose $g\in\Gamma$ with $\mu(g)>0$. 
By symmetry, $\mu(g)=\mu(g^{-1})>0$.
Hence, the event $\mathcal{A}$ that the root has $k$ children, two of them move to $g$, each of these two particles has $k$ children, and one child of each moves back to $e$ has probability $\varrho>0$. 
On $\mathcal{A}$, there are two particles at $e$ whose descendant BRWs are independent copies of the original BRW started at $e$.

Let $L_1$ and $L_2$ denote the local times of these two descendant BRWs. 
On the event $\mathcal{A}$,  we have $L(z)\geq L_1(z)+L_2(z)$ for all $z\in\Gamma$. 
Therefore, for all $s>1$,
\begin{align*}
	\left(s^{L(x)}-1\right) \left(s^{L(y)}-1\right) &\geq \left(s^{L_1(x)+L_2(x)}-1\right) \left(s^{L_1(y)+L_2(y)}-1\right)  \\
	&\geq \left(s^{L_1(x)}-1\right) \left(s^{L_1(y)}-1\right) \left(s^{L_2(x)}-1\right)\left (s^{L_2(y)}-1\right).
\end{align*}
By the branching property and the conditional independence of the two descendant BRWs, we obtain
\begin{equation*}
	 F(s)\geq \varrho F(s)^2.
\end{equation*}
Since $F(s)<\infty$ for $ 1<s<e^{c^*}$, it follows that $F(s)\leq \varrho^{-1}$.
Letting $s\uparrow e^{c^*}$ and using the monotone convergence theorem, we obtain
\begin{equation*}
	 F(e^{c^*})\leq\varrho^{-1}<\infty.
\end{equation*}
On the other hand, by \eqref{eq:joint-tail-generating-function} and the lower bound in \eqref{eq:equal-exponential-rates}, we have 
\begin{equation*}
	 F(e^{c^*}) \geq  c e^{-2c^*} (e^{c^*}-1)^2 \Phi_r(e;x,y) \sum_{n=1}^{\infty} \sum_{m=1}^{\infty} 1 =\infty,
\end{equation*}
which is a contradiction. 
Therefore, $c_2> C_2$.
\end{proof}

Theorem~\ref{thm:two-point-tail-estimate} now follows from Propositions~\ref{prop:two-sided-estimates} and~\ref{prop:strict-exponential-constants}.

\subsection{The subcritical case on finitely generated groups}  \label{subsec:subcritical-green-convolution}
We now assume $r\in(0,R)$. 
Thus, $G_r(x,y)<\infty$ for all $x,y\in \Gamma$.
Recall that $\rho=R^{-1}$ is the spectral radius of the underlying RW, so that
\begin{equation*}
    r\rho<1.
\end{equation*}
By the left invariance of the Green function, we have the following triangle estimate.

\begin{lemma} [{\cite[Lemma~4.3]{GL26}}]  \label{lem:subcritical-triangle-green}
For every $r\in(0,R)$, there exists a constant $C=C(r,\mu)<\infty$ such that, for all $a,x\in\Gamma$,
\begin{equation} \label{eq:subcritical-triangle-green}
    \sum_{z\in\Gamma}G_r(a,z)^2 G_r(z,x)  \leq C G_r(a,x).
\end{equation}
\end{lemma}

Lemma~\ref{lem:subcritical-triangle-green} verifies the condition~\eqref{eq:criterion-triangle}.
It remains to verify the condition~\eqref{eq:criterion-mixed}.

\begin{proposition} \label{prop:subcritical-G-Phi-convolution}
For every $r\in(0,R)$, there exists a constant $C=C(r,\mu)<\infty$ such that, for all $a,x,y\in\Gamma$,
\begin{equation} \label{eq:subcritical-G-Phi-convolution}
    \sum_{z\in\Gamma}G_r(a,z)G_r(z,x)\Phi_r(z;x,y) \leq C\Phi_r(a;x,y).
\end{equation}
\end{proposition}

\begin{proof}
Expanding the left-hand side of \eqref{eq:subcritical-G-Phi-convolution}, we obtain
\begin{align}
    &\sum_{z\in\Gamma}G_r(a,z)G_r(z,x)\Phi_r(z;x,y) \notag\\
    &\qquad=\sum_{z,u\in\Gamma}\sum_{i,j,k,l,s\geq0} r^{i+j+k+l+s}  p_i(a,z)p_j(z,x)p_k(z,u)p_l(u,x)p_s(u,y).  \label{eq:subcritical-G-Phi-expanded}
\end{align}
The indices $j,k,l$ correspond to the three Green functions forming the triangle with vertices $z,x,u$. 
Since at least one of $j,k,l$ is maximal, it suffices to consider the following three cases.

\smallskip
\noindent\textbf{Case 1. $j\geq k\vee l$.}
Let $I_1$ denote the contribution to \eqref{eq:subcritical-G-Phi-expanded} from the indices satisfying $j\geq k\vee l$.
By $p_j(z,x)\leq\rho^j$ and the Chapman--Kolmogorov identity, we obtain
\begin{align*}
    I_1 &\leq \sum_{\substack{i,j,k, l,s\geq0\\j\geq k\vee l}} (r\rho)^j r^{i+k+l+s}  \sum_{u\in\Gamma}p_{i+k}(a,u)p_l(u,x)p_s(u,y)\\
     &\leq C\sum_{i,k,l,s\geq0} (r\rho)^{k\vee l}r^{i+k+l+s} \sum_{u\in\Gamma}p_{i+k}(a,u)p_l(u,x)p_s(u,y),
\end{align*}
where we used $ \sum_{j\geq k\vee l}(r\rho)^j \leq C(r\rho)^{k\vee l}$ for some constant $C=C(r\rho)<\infty$, since $r\rho<1$.
Setting $n=i+k$ and observing that
\begin{equation*}
    \sum_{k=0}^{n}(r\rho)^{k\vee l} \leq (l+1)(r\rho)^l+\sum_{k=l+1}^{\infty}(r\rho)^k \leq C,
\end{equation*}
for all $n,l \geq 0$, where $C=C(r\rho)<\infty$, we obtain
\begin{align*}
    I_1 &\leq C \sum_{n,l,s\geq0 } r^{n+l+s} \sum_{u\in\Gamma} p_{n}(a,u) p_l(u,x)  p_s(u,y) \sum_{k=0}^n(r\rho)^{k\vee l}  \\
        &\leq  C  \Phi_r (a;x,y).
\end{align*}

\smallskip 
\noindent\textbf{Case 2. $k\geq j\vee l$.}
Let $I_2$ denote the corresponding contribution. Using $p_k(z,u)\leq\rho^k$, the symmetry of the transition probabilities, and the Chapman--Kolmogorov identity, we obtain
\begin{align*}
    I_2 &\leq \sum_{\substack{i,j,k,l,s\geq0\\k\geq j\vee l}} (r\rho)^k r^{i+j+l+s}p_{i+j}(a,x)p_{l+s}(x,y)\\
        &\leq C\sum_{i,j,l,s\geq0}  (r\rho)^{j\vee l}r^{i+j+l+s}p_{i+j}(a,x)p_{l+s}(x,y).
\end{align*}
Since $(r\rho)^{j\vee l} \leq (r\rho)^{(j+l)/2}$,
\begin{align*}
    I_2 &\leq C \left(\sum_{i,j\geq0}(r\rho)^{j/2}r^{i+j}p_{i+j}(a,x)\right) \left(\sum_{l,s\geq0}(r\rho)^{l/2}r^{l+s}p_{l+s}(x,y)\right)\\
        &\leq C G_r(a,x)G_r(x,y),
\end{align*}
where the last inequality follows from the summability of the geometric series $\sum_{j\geq0}(r\rho)^{j/2}$.
On the other hand, considering only the term $z=x$ in the definition of $\Phi_r$, we have $\Phi_r(a;x,y)\geq G_r(e,e)G_r(a,x)G_r(x,y)$. Hence $I_2\leq C\Phi_r(a;x,y)$.

\smallskip \noindent\textbf{Case 3. $l\geq j\vee k$.}
Let $I_3$ denote the corresponding contribution. 
Using $p_l(u,x)\leq\rho^l$ and the Chapman--Kolmogorov identity, we obtain
\begin{align*}
    I_3 &\leq \sum_{\substack{i,j,k,l,s\geq0\\ l\geq j\vee k}} (r\rho)^l r^{i+j+k+s} \sum_{z\in\Gamma}p_i(a,z)p_j(z,x)p_{k+s}(z,y)\\
        &\leq C\sum_{i,j,k,s\geq0} (r\rho)^{j\vee k}r^{i+j+k+s} \sum_{z\in\Gamma}p_i(a,z)p_j(z,x)p_{k+s}(z,y).
\end{align*}
Setting $n=k+s$ and arguing as in Case~1, we obtain $I_3  \leq C\Phi_r(a;x,y)$.

Combining the three cases, we get the desired result.
\end{proof}

\begin{proof}[Proof of Corollary~\ref{cor:subcritical-two-point-tail-estimate}]
Lemma~\ref{lem:subcritical-triangle-green} and Proposition~\ref{prop:subcritical-G-Phi-convolution} verify the conditions~\eqref{eq:criterion-triangle} and \eqref{eq:criterion-mixed}, respectively.
The result therefore follows from Theorem~\ref{thm:two-point-tail-estimate}. 
\end{proof}

\section{The hyperbolic critical case} \label{sec:hyperbolic-case}
Throughout this section, let $\Gamma$ be a non-elementary hyperbolic group, and let $\mu$ be an admissible symmetric probability measure on $\Gamma$ with superexponential tails. 
Suppose that the offspring distribution $\nu$ satisfies Assumption~\ref{ass:offspring} with mean $r=R$.

For $r=R$, the inequality $r\rho<1$ is no longer available.
We instead use Ancona inequalities together with geodesic projection estimates from hyperbolic geometry to verify the conditions~\eqref{eq:criterion-triangle} and~\eqref{eq:criterion-mixed} at $r=R$.

The corresponding one-point local time tail estimate was established in \cite[Theorem~1.2]{GL26}: there exist positive constants $c_1,c_2,C_1,C_2$, depending only on $\nu$ and $\mu$, such that, for all $a,x\in\Gamma$ and all integers $n\geq1$,
\begin{equation} \label{eq:critical-one-point-tail-hyperbolic}
    c_1e^{-c_2n}G_R(a,x) \leq \pr_a\bigl(L(x)\geq n\bigr) \leq C_1e^{-C_2n}G_R(a,x).
\end{equation}

\subsection{Hyperbolic groups}
Let $d$ be the word metric associated with a finite symmetric generating set. 
We write $|x|:=d(e,x)$ for $x\in\Gamma$.
For $x, y \in \Gamma$, let $[x,y]$ be an arbitrarily chosen geodesic segment connecting $x$ and $y$. 

The Gromov product of $x$ and $y$ with respect to $z$ is defined by 
\begin{equation}
	(x|y)_z :=\frac{1}{2}[d(x,z)+d(y,z)-d(x,y)].
\end{equation}
For $\delta \geq 0$, we say that $(\Gamma, d)$ is $\delta$-\textbf{hyperbolic} if, for all $x, y, z, w \in \Gamma$, 
\begin{equation}
	(x|y)_w \geq \min\{ (x|z)_w,  (z|y)_w \} -\delta.
\end{equation}
Say that $(\Gamma, d)$ is hyperbolic if it is $\delta$-hyperbolic for some $\delta \geq 0$.
From now on, we always assume that $\Gamma$ is \textbf{non-elementary}, that is $\Gamma$ is  neither finite nor virtually $\mathbb{Z}$.

If $(\Gamma, d)$ is $\delta$-hyperbolic, then every geodesic triangle is $4\delta$-thin, that is, each side of the triangle is contained in the $4\delta$-neighborhood of the union of the other two sides. 

Let $\kappa \geq 0$.
A nonempty subset $Q\subset \Gamma$ is called \textbf{$\kappa$-quasiconvex} if every geodesic segment connecting two points of $Q$ is contained in the $\kappa$-neighborhood $N_{\kappa}(Q):=\{z\in\Gamma:d(z,Q)\leq\kappa\}$ of $Q$.
For a quasiconvex subset $Q\subset \Gamma$ and $u\in \Gamma$, let $\pi_Q(u)$ be a \textbf{nearest point projection} of $u$ onto $Q$, that is, $d(u, \pi_Q(u))= d(u, Q)$. 
If there is more than one such point, we choose one arbitrarily. 

Quasiconvex subsets satisfy the following coarse projection property, which follows from \cite[Chapitre 10, Proposition 2.1]{CDP90}.

\begin{proposition} [{\cite[Chapitre 10, Proposition 2.1]{CDP90}}]  \label{prop:projection-quasiconvex-group}
Let $Q\subset \Gamma$ be a nonempty $\kappa$-quasiconvex subset. 
For $u\in \Gamma$, let $\pi_Q(u)$ be a nearest point projection of $u$ onto $Q$.  
Then there exists a constant $D=D(\delta,\kappa)<\infty$ such that, for all $u\in\Gamma$ and $w\in Q$, 
\begin{equation*}
    d \left( \pi_Q(u),[u,w] \right)\leq D.
\end{equation*}
In particular, every geodesic segment from $u$ to a point $w\in Q$ passes through a uniformly bounded neighborhood of $\pi_Q(u)$.
\end{proposition}

We first note that every geodesic segment is $4\delta$-quasiconvex. 
Indeed, let $a,b\in[x,y]$. 
The subsegment of $[x,y]$ connecting $a$ and $b$ is itself a geodesic segment. 
Viewing this subsegment and any other geodesic segment connecting $a$ and $b$ as two sides of a degenerate geodesic triangle, the $4\delta$-thinness of geodesic triangles implies that every geodesic segment connecting $a$ and $b$ is contained in the $4\delta$-neighborhood of this subsegment.

Furthermore, for $x,y,z\in\Gamma$, let $\pi(z):=\pi_{[x,y]}(z)$ and set $ T:=[z,\pi(z)]\cup[\pi(z),x]\cup[\pi(z),y]$, where $[\pi(z),x]$ and $[\pi(z),y]$ are the corresponding subsegments of the chosen geodesic $[x,y]$.
Then $T$ is also $4\delta$-quasiconvex.  
Indeed, $T$ is the union of three geodesic segments with the common endpoint $\pi(z)$. 
Let $a,b\in T$. 
If $a$ and $b$ lie on the same geodesic segment, then the preceding argument shows that every geodesic segment connecting $a$ and $b$ is contained in the $4\delta$-neighborhood of $T$. 
If $a$ and $b$ lie on two different geodesic segments, then $[a,\pi(z)]\cup[\pi(z),b]\subseteq T$. 
By the $4\delta$-thinness of geodesic triangles, every geodesic segment connecting $a$ and $b$ is contained in the $4\delta$-neighborhood of $[a,\pi(z)]\cup[\pi(z),b]$, and hence in the $4\delta$-neighborhood of $T$. 
Therefore, Proposition~\ref{prop:projection-quasiconvex-group} implies the following corollary.

\begin{corollary} \label{cor:geodesic-projection}
For $x,y,z\in\Gamma$, let $\pi(z):=\pi_{[x,y]}(z)$ be a nearest point projection of $z$ onto the geodesic segment $[x,y]$. Then there exists a constant $D=D(\delta)<\infty$ such that
\begin{equation*}
    d\bigl(\pi(z),[z,x]\bigr)\leq D,
    \qquad
    d\bigl(\pi(z),[z,y]\bigr)\leq D.
\end{equation*}
Moreover, let $T:=[z,\pi(z)]\cup[\pi(z),x]\cup[\pi(z),y]$ and let $\pi_T(u)$ be a nearest point projection of $u$ onto $T$. Then for every $u\in\Gamma$ and $w\in T$,
\begin{equation*}
    d\bigl(\pi_T(u),[u,w]\bigr)\leq D.
\end{equation*}
\end{corollary}

\subsection{Critical Green function estimates}  \label{subsec:critical-Green-convolution}
Recall that $R$ is the radius of convergence of the Green function $G_r(x, y)$. 
Since $\Gamma$ is non-elementary, it is non-amenable, and hence by \cite[Corollary 12.5(a)]{Woe00}, we know that $R>1$. 
Moreover, by \cite[Theorem~7.8]{Woe00}, for all $x,y\in\Gamma$, 
\begin{equation}
    G_R(x,y)<\infty.
\end{equation}

We first record the following  comparison of Green function at $R$.

\begin{lemma} \label{lem:critical-Ancona-inequality}
For every $D>0$, there exists a constant $C=C(D,R,\mu)<\infty$ such that, for all $x,y,z\in\Gamma$ satisfying $d\bigl(y,[x,z]\bigr)\leq D$,
\begin{equation} \label{eq:critical-Ancona-inequality}
    C^{-1} G_R(x,y)G_R(y,z) \leq G_R(x,z)\leq C G_R(x,y)G_R(y,z).
\end{equation}	
\end{lemma}

The upper bound follows from the Ancona inequality at $R$ \cite[Theorem~1.3]{Gou15}, together with Harnack inequalities.
The lower bound follows from $G_R(x,y)G_R(y,z)\leq G_R(e,e)G_R(x,z)$, which holds for all $x,y,z \in \Gamma$. 

We next recall the exponential decay of the Green function at $R$.

\begin{lemma} [{\cite[Proposition~5.3]{GL26}}] \label{lem:critical-Green-exponential-decay}
There exist constants $C=C(R, \mu)<\infty$ and $\kappa=\kappa(R, \mu)>0$ such that, for all $x,y\in\Gamma$,
\begin{equation*}  
    G_R(x,y)\leq Ce^{-\kappa d(x,y)}.
\end{equation*}
In particular, $G_R$ is uniformly bounded on $\Gamma\times\Gamma$.
\end{lemma}

We will also use the following two summability estimates from \cite[Lemmas~5.5, 5.6]{GL26}.

\begin{lemma} \label{lem:summability-estimates}
There exists a constant $C=C(R, \mu)<\infty$ such that the following hold:
\begin{itemize}
	\item[(i)]  $\sum_{z\in\Gamma} G_R(e,z)^3 \leq C$;
	\item[(ii)] for all $a,x\in \Gamma$, $ \sum_{z\in\Gamma} G_R(a,z)^2 G_R(z,x)\leq C G_R(a,x)$.
\end{itemize}	
\end{lemma}

Lemma~\ref{lem:summability-estimates} (ii) verifies the triangle condition~\eqref{eq:criterion-triangle}.
It remains to verify the mixed convolution condition~\eqref{eq:criterion-mixed}. 
We first establish a projection estimate that identifies the spatial factor $\Phi_R(a;x,y)$.

\begin{proposition} \label{prop:critical-Phi-projection}
There exists a constant $C=C(R,\mu)<\infty$ such that, for all $a,x,y\in\Gamma$,
\begin{equation} \label{eq:critical-Phi-projection}
    \Phi_R(a;x,y) \leq C G_R(a,\pi(a))G_R(\pi(a),x)G_R(\pi(a),y),
\end{equation}
where $\pi(a):=\pi_{[x,y]}(a)$ is a nearest point projection of $a$ onto the geodesic segment $[x,y]$. 
In particular, $ \Phi_R(a;x,y)<\infty$.
\end{proposition}

\begin{proof}
Write $\pi:=\pi(a)$ and set $ T:=[a,\pi]\cup[\pi,x]\cup[\pi,y]$.
Fix $z\in\Gamma$, and let $v:=\pi_T(z)$ be a nearest point projection of $z$ onto $T$.
By Corollary~\ref{cor:geodesic-projection}, there exists a constant $D<\infty$ such that
\begin{equation*}
    d(v,[z,a])\leq D,\qquad d(v,[z,x])\leq D,\qquad \text{and} \qquad d(v,[z,y])\leq D.
\end{equation*}
Hence, by the upper bound in Lemma~\ref{lem:critical-Ancona-inequality} and the symmetry of $G_R$,
\begin{equation} \label{eq:critical-Phi-projection-first-step}
    G_R(a,z)G_R(z,x)G_R(z,y)\leq C G_R(a,v)G_R(v,x)G_R(v,y)G_R(v,z)^3.
\end{equation}

We next compare $G_R(a,v)G_R(v,x)G_R(v,y)$ with the corresponding product at $\pi$. 
Since $ v\in T$, we distinguish three cases; see Figure~\ref{fig:critical-Phi-projection-cases}.

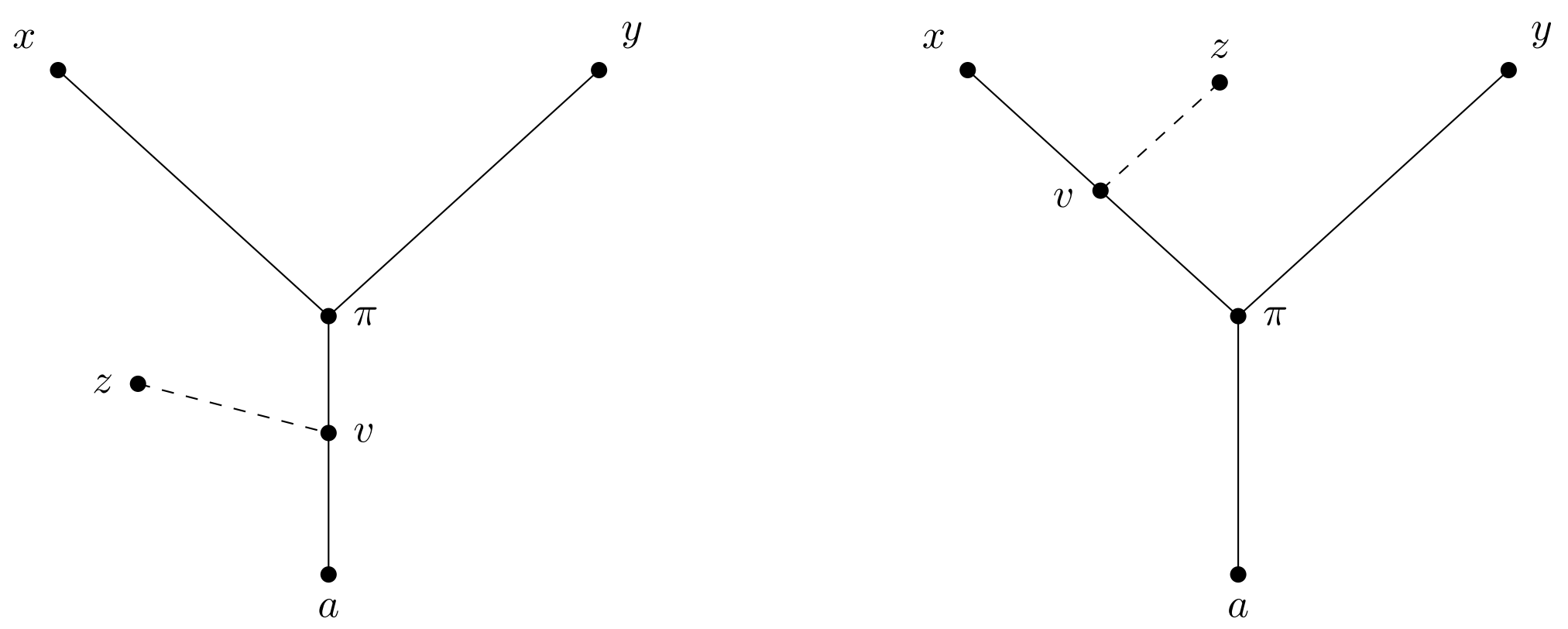
\begin{figure}[htbp]
\centering
\begin{tikzpicture}[
scale=1.05,
every node/.style={font=\small},
dot/.style={circle,fill,inner sep=1.4pt}
]

\begin{scope}[xshift=0cm]

\coordinate (p) at (0,0);
\coordinate (x) at (-2.2,2.0);
\coordinate (y) at (2.2,2.0);
\coordinate (a) at (0,-2.1);

\coordinate (v) at (0,-0.95);
\coordinate (z) at (-1.55,-0.55);

\draw (x)--(p)--(y);
\draw (a)--(p);
\draw[dashed] (z)--(v);

\node[dot,label=above left:$x$] at (x) {};
\node[dot,label=above right:$y$] at (y) {};
\node[dot,label=below:$a$] at (a) {};
\node[dot,label=right:$\pi$] at (p) {};
\node[dot,label=right:$v$] at (v) {};
\node[dot,label=left:$z$] at (z) {};

\node at (0,-2.8) {(a) \(v\in[a,\pi]\)};

\end{scope}

\begin{scope}[xshift=7.4cm]

\coordinate (p) at (0,0);
\coordinate (x) at (-2.2,2.0);
\coordinate (y) at (2.2,2.0);
\coordinate (a) at (0,-2.1);

\coordinate (v) at (-1.12,1.02);
\coordinate (z) at (-0.15,1.9);

\draw (x)--(p)--(y);
\draw (a)--(p);
\draw[dashed] (z)--(v);

\node[dot,label=above left:$x$] at (x) {};
\node[dot,label=above right:$y$] at (y) {};
\node[dot,label=below:$a$] at (a) {};
\node[dot,label=right:$\pi$] at (p) {};
\node[dot] at (v) {};
\node[fill=white,inner sep=1pt] at (-1.42,0.96) {$v$};
\node[dot,label=above:$z$] at (z) {};

\node at (0,-2.8) {(b) \(v\in[\pi,x]\)};

\end{scope}

\end{tikzpicture}
\caption{Two possible locations of $v=\pi_T(z)$. The case $v\in[\pi,y]$ is symmetric to~(b) after interchanging $x$ and $y$.}
\label{fig:critical-Phi-projection-cases}
\end{figure}

\smallskip
\noindent\textbf{Case 1. $v\in[a,\pi]$.}
Since $\pi$ is a nearest point projection of $a$ onto $[x,y]$, it is also a nearest point projection of $v$ onto $[x,y]$.
Hence, by Corollary~\ref{cor:geodesic-projection}, we have $ d(\pi,[v,x])\leq D$ and $d(\pi,[v,y])\leq D$. 
Therefore, the upper bound in Lemma~\ref{lem:critical-Ancona-inequality} gives
\begin{equation*}
    G_R(v,x)\leq C G_R(v,\pi)G_R(\pi,x) \qquad \text{and} \qquad  G_R(v,y)\leq C G_R(v,\pi)G_R(\pi,y).
\end{equation*}
Using the lower bound in Lemma~\ref{lem:critical-Ancona-inequality}, we then obtain
\begin{align*}
    G_R(a,v)G_R(v,x)G_R(v,y)
    &\leq C G_R(a,v)G_R(v,\pi)^2G_R(\pi,x)G_R(\pi,y)\\
    &\leq C G_R(a,\pi)G_R(\pi,x)G_R(\pi,y)G_R(v,\pi).
\end{align*}

\smallskip
\noindent\textbf{Case 2. $v\in[\pi,x]$.}
Let $[v,y]$ be the subsegment of $[x,y]$ connecting $v$ and $y$. 
Then $\pi$ is also a nearest point projection of $a$ onto the subsegment $[v,y]$. 
Hence, by Corollary~\ref{cor:geodesic-projection}, we have $d(\pi,[a,v])\leq D$.
Therefore, by the upper bound in Lemma~\ref{lem:critical-Ancona-inequality}, 
\begin{equation*}
    G_R(a,v)\leq C G_R(a,\pi)G_R(\pi,v) \qquad \text{and} \qquad G_R(v,y)\leq C G_R(v,\pi)G_R(\pi,y).
\end{equation*}
Using the lower bound in Lemma~\ref{lem:critical-Ancona-inequality} again, we obtain 
\begin{align*}
    G_R(a,v)G_R(v,x)G_R(v,y)
    &\leq C G_R(a,\pi)G_R(\pi,v)G_R(v,x)G_R(v,\pi)G_R(\pi,y)\\
    &\leq C G_R(a,\pi)G_R(\pi,x)G_R(\pi,y)G_R(v,\pi).
\end{align*}

\smallskip
\noindent\textbf{Case 3. $v\in[\pi,y]$.}
This case is symmetric to Case~2 after interchanging $x$ and $y$.

Combining the three cases, we have, for every $v\in T$,
\begin{equation} \label{eq:critical-Phi-projection-product}
    G_R(a,v)G_R(v,x)G_R(v,y)
    \leq C G_R(a,\pi)G_R(\pi,x)G_R(\pi,y)G_R(v,\pi).
\end{equation}
Using \eqref{eq:critical-Phi-projection-first-step} and \eqref{eq:critical-Phi-projection-product}, and decomposing the sum over $z\in\Gamma$ according to $v=\pi_T(z)$, we obtain
\begin{equation} \label{eq:Phi-R}
	  \Phi_R(a;x,y) \leq C G_R(a,\pi)G_R(\pi,x)G_R(\pi,y) \sum_{v\in T}G_R(v,\pi) \sum_{\substack{z\in\Gamma \\ \pi_T(z)=v}}G_R(v,z)^3.
\end{equation}
By left invariance and Lemma~\ref{lem:summability-estimates} (i),
\begin{equation} \label{eq:cubic-sum}
    \sum_{\substack{z\in\Gamma\\\pi_T(z)=v}}G_R(v,z)^3 \leq \sum_{z\in\Gamma}G_R(e,z)^3 \leq C.
\end{equation}
Moreover, since $T$ is the union of three geodesic segments,  Lemma~\ref{lem:critical-Green-exponential-decay} gives
\begin{equation} \label{eq:geodesic-sum}
    \sum_{v\in T} G_R(v,\pi) \leq 3C\sum_{k=0}^{\infty}e^{-\kappa k} <\infty.
\end{equation}
Therefore, combining \eqref{eq:Phi-R}, \eqref{eq:cubic-sum}, and \eqref{eq:geodesic-sum}, we obtain the desired bound and conclude that $\Phi_R(a;x,y)<\infty$. 
\end{proof}

The upper bound in Proposition~\ref{prop:critical-Phi-projection} is complemented by an immediate lower bound.
Taking $z=\pi(a)$ in the definition of $\Phi_R(a;x,y)$, we have
\begin{equation*}
	 \Phi_R(a;x,y) \geq  G_R(a,\pi(a))G_R(\pi(a),x)G_R(\pi(a),y).
\end{equation*}
Thus combining this with Proposition~\ref{prop:critical-Phi-projection} yields the following corollary.

\begin{corollary} \label{cor:critical-Phi-projection-asymptotic}
For all $a,x,y \in \Gamma$,
\begin{equation} \label{eq:critical-Phi-projection-asymptotic}
	 \Phi_R(a;x,y) \asymp  G_R(a,\pi(a))G_R(\pi(a),x)G_R(\pi(a),y),
\end{equation}	
where the implicit constants depend only on $R$ and $\mu$.
\end{corollary}

We now verify the mixed convolution condition~\eqref{eq:criterion-mixed}.

\begin{proposition} \label{prop:critical-G-Phi-convolution}
There exists a constant $C=C(R,\mu)<\infty$ such that, for all $a,x,y\in\Gamma$,
\begin{equation} \label{eq:critical-G-Phi-convolution-x}
     \sum_{z\in\Gamma} G_R(a,z)G_R(z,x)\Phi_R(z;x,y)  \leq C\Phi_R(a;x,y).
\end{equation}
\end{proposition}

\begin{proof}
For $z\in\Gamma$, let $\pi(z):=\pi_{[x,y]}(z)$. 
By Corollary~\ref{cor:geodesic-projection} and the upper bound in Lemma~\ref{lem:critical-Ancona-inequality},
\begin{equation*}
    G_R(z,x)\leq C G_R(z,\pi(z)) G_R(\pi(z),x).
\end{equation*}
Proposition~\ref{prop:critical-Phi-projection} also gives
\begin{equation*}
    \Phi_R(z;x,y) \leq C G_R(z,\pi(z))G_R(\pi(z),x)G_R(\pi(z),y).
\end{equation*}
Therefore, decomposing the sum according to $u=\pi(z)\in[x,y]$, and then using the symmetry of $G_R$ together with Lemma~\ref{lem:summability-estimates} (ii), we obtain
\begin{align*}
    &\sum_{z\in\Gamma}G_R(a,z)G_R(z,x)\Phi_R(z;x,y)
      \leq C \sum_{u\in[x,y]} G_R(u,x)^2G_R(u,y) \sum_{\substack{z\in\Gamma\\\pi(z)=u}} G_R(a,z)G_R(z,u)^2 \\
    & \qquad \leq C \sum_{u\in[x,y]} G_R(a,u) G_R(u,x)^2 G_R(u,y) \leq C \Phi_R(a;x,y),
\end{align*} 
where the last inequality follows from the uniform boundedness of $G_R$.
This completes the proof.
\end{proof}

\subsection{Critical two-point local time tail estimates} \label{subsec:critical-two-point-tail-estimate}

\begin{proof}[Proof of Corollary~\ref{cor:critical-two-point-tail-estimate-hyperbolic}]
Lemma~\ref{lem:summability-estimates}(ii) and Proposition~\ref{prop:critical-G-Phi-convolution} verify the conditions~\eqref{eq:criterion-triangle} and~\eqref{eq:criterion-mixed}, respectively. 
Since $G_R(x,y)<\infty$ for all $x,y\in\Gamma$, Theorem~\ref{thm:two-point-tail-estimate} yields the desired two-point tail estimate. 
Moreover, Corollary~\ref{cor:critical-Phi-projection-asymptotic} gives
\begin{equation*}
	\Phi_R(a;x,y) \asymp  G_R(a,\pi(a))G_R(\pi(a),x)G_R(\pi(a),y).
\end{equation*} 
This completes the proof.
\end{proof}

By Corollary~\ref{cor:geodesic-projection} and Lemma~\ref{lem:critical-Ancona-inequality}, we have
\begin{equation*}
	G_R(a,x)\asymp G_R(a,\pi(a))G_R(\pi(a),x) \qquad \text{and} \qquad G_R(a,y)\asymp G_R(a,\pi(a))G_R(\pi(a),y).
\end{equation*}	
Combining these estimates with Corollaries~\ref{cor:critical-Phi-projection-asymptotic} and~\ref{cor:critical-two-point-tail-estimate-hyperbolic} yields the following corollary.

\begin{corollary} \label{cor:critical-two-point-tail-geometric}
There exist positive constants $c_1,c_2,C_1,C_2$ depending only on $\nu$ and $\mu$ such that, for all $a,x,y\in\Gamma$ with $x\neq y$ and all integers $n,m\geq 1$,
\begin{equation*}
    c_1e^{-c_2(n+m)} \frac{G_R(a,x)G_R(a,y)}{G_R(a,\pi)} \leq \pr_a \bigl(L(x)\geq n, L(y)\geq m\bigr)  \leq C_1e^{-C_2(n+m)}  \frac{G_R(a,x)G_R(a,y)}{G_R(a,\pi)},
\end{equation*}
where $\pi:=\pi_{[x,y]}(a)$.
Moreover, any such choice of constants necessarily satisfies $c_2>C_2$.
\end{corollary}

Taking $n=m=1$ in  Corollary~\ref{cor:critical-two-point-tail-geometric} and using the one-point estimate \eqref{eq:critical-one-point-tail-hyperbolic}, we obtain
\begin{equation*}
    \frac{\pr_a(L(x)\geq1,L(y)\geq1)}  {\pr_a(L(x)\geq1)\pr_a(L(y)\geq1)} \asymp  \frac{1}{G_R(a,\pi)}.
\end{equation*}
Consequently, $G_R(a,\pi)^{-1}$ quantifies the spatial dependence between the two visiting events up to multiplicative constants. 
By Lemma~\ref{lem:critical-Green-exponential-decay}, this factor is at least $ce^{\kappa d(a,\pi)}$ for suitable constants $c,\kappa>0$.

\subsection{An application to the growth of visited sites on spheres} \label{subsec:visited-sphere-second-moment}

We conclude this section with an application of the two-point estimate to the growth of the number of visited sites on spheres, recovering the corresponding result in \cite[Theorem~4.1]{SWX23}. 
Let $\mathbb{S}_n:=\{x\in\Gamma:d(x,e)=n\}$ be the sphere of radius $n$ centered at $e$, and let
\begin{equation*}
    M_n:=\sum_{x\in\mathbb{S}_n}\mathbf{1}_{\{L(x)\geq1\}}
\end{equation*}
denote the number of sites in $\mathbb{S}_n$ ever visited by the BRW.
Set
\begin{equation*}
	 H_n(r):=\sum_{x\in\mathbb{S}_n}G_r(e,x)\quad \text{ and } \quad
     H(r):=\limsup_{n\to\infty}H_n(r)^{1/n}.
\end{equation*}
By \cite[Theorem~3.1]{SWX23}, $H(R)>1$, and there exists a constant $C>1$ such that, for all integers $n \geq 0$,
\begin{equation} \label{eq:critical-Hn-growth}
    C^{-1}H(R)^n\leq H_n(R)\leq CH(R)^n. 
\end{equation}

We first establish a summability estimate for the spatial factor in Corollary~\ref{cor:critical-two-point-tail-geometric}. 
For $x\neq y$, write $ \pi_{x,y}:=\pi_{[x,y]}(e)$.

\begin{lemma} \label{lem:critical-sphere-correlation-sum}
There exists a constant $C>0$ such that, for all integers $n\geq1$,
\begin{equation*}
    \sum_{\substack{x,y\in\mathbb{S}_n\\x\neq y}} \frac{G_R(e,x)G_R(e,y)}{G_R(e,\pi_{x,y})} \leq CH_n(R)^2.
\end{equation*}
\end{lemma}

\begin{proof}
Fix an integer $A\geq1$ larger than the constant in Corollary~\ref{cor:geodesic-projection}, and for $u\in\Gamma$ define
\begin{equation*}
    S_{n,A}(u):=\sum_{\substack{x\in\mathbb{S}_n\\d(u,[e,x])\leq A}}G_R(u,x).
\end{equation*}
For $x\in\mathbb{S}_n$, if $d(u,[e,x])\leq A$, then $\left|d(u,x)-(n-|u|)_+\right|\leq 2A$.
Hence, by the left invariance of the Green function and \eqref{eq:critical-Hn-growth},
\begin{equation} \label{eq:critical-SnA-upper}
    S_{n,A}(u) \leq \sum_{\substack{k\geq0 \\ |k-(n-|u|)_+|\leq2A}}H_k(R) \leq C_A H(R)^{(n-|u|)_+}.    
\end{equation}

For $x\neq y$, Corollary~\ref{cor:geodesic-projection} and Lemma~\ref{lem:critical-Ancona-inequality} yield
\begin{equation*}
    \frac{G_R(e,x)G_R(e,y)}{G_R(e,\pi_{x,y})}\asymp G_R(e,\pi_{x,y})G_R(\pi_{x,y},x)G_R(\pi_{x,y},y).
\end{equation*}
Note that $\pi_{x,y}$ lies within distance $A$ of both $[e,x]$ and $[e,y]$.  
Therefore, decomposing the sum according to $u=\pi_{x,y}$, we obtain
\begin{equation*}
       \sum_{\substack{x,y\in\mathbb{S}_n\\x\neq y}} \frac{G_R(e,x)G_R(e,y)}{G_R(e, \pi_{x,y})} \leq C\sum_{u\in\Gamma}G_R(e,u)S_{n,A}(u)^2.
\end{equation*}
A nonzero term on the right-hand side necessarily satisfies $|u|\leq n+A$. 
Hence, by \eqref{eq:critical-SnA-upper} and \eqref{eq:critical-Hn-growth},
\begin{align*}
   \sum_{u\in\Gamma}G_R(e,u)S_{n,A}(u)^2
    &\leq C_A\sum_{k=0}^{n+A}H(R)^kH(R)^{2(n-k)_+}\\
    &\leq C_AH(R)^{2n} \leq C_AH_n(R)^2,
\end{align*}	
where we used that $H(R)>1$. This proves the lemma.
\end{proof}

\begin{corollary}  \label{cor:critical-visited-sphere-growth}
On the survival event,
\begin{equation} \label{eq:M-as}
    \limsup_{n\to\infty}M_n^{1/n}=H(R)  
\end{equation}
almost surely. Consequently,
\begin{equation*}
    \limsup_{n\to\infty}\frac{1}{n}\log M_n  = \limsup_{n\to\infty}\frac{1}{n}\log\E[M_n] = \log H(R),
\end{equation*}
where we use the convention $\log 0=-\infty$.
\end{corollary}

\begin{proof}
By \eqref{eq:critical-Hn-growth} and the one-point estimate \eqref{eq:critical-one-point-tail-hyperbolic},
\begin{equation} \label{eq:critical-visited-sphere-first-moment-growth}
     \E[M_n] \asymp H_n(R) \asymp H(R)^n.
\end{equation}
Thus, for every $\varepsilon>0$,  by Markov's inequality, we obtain
\begin{equation*}
    \pr\bigl(M_n^{1/n}\geq H(R)+\varepsilon\bigr)  \leq C\left(\frac{H(R)}{H(R)+\varepsilon}\right)^n.
\end{equation*}
Consequently, the Borel--Cantelli lemma yields
\begin{equation*}
    \limsup_{n\to\infty}M_n^{1/n}\leq H(R) \qquad\text{a.s.}
\end{equation*}

For the reverse inequality, by \eqref{eq:critical-one-point-tail-hyperbolic}, Corollary~\ref{cor:critical-two-point-tail-geometric}, Lemma~\ref{lem:critical-sphere-correlation-sum} and \eqref{eq:critical-visited-sphere-first-moment-growth}, we have 
\begin{align*}
    \E[M_n^2]
    &=\sum_{x\in\mathbb{S}_n}\pr\bigl(L(x)\geq1\bigr) +\sum_{\substack{x,y\in\mathbb{S}_n\\x\neq y}} \pr\bigl(L(x)\geq1,L(y)\geq1\bigr)\\
    &\leq CH_n(R)+CH_n(R)^2
    \leq C\E[M_n]^2.
\end{align*}
Hence, the Paley--Zygmund inequality yields
\begin{equation*}
    \pr\left(M_n\geq\frac{1}{2}\E[M_n]\right)\geq c>0,
\end{equation*}
uniformly in $n$.

Fix $0<\varepsilon<H(R)-1$.  
By \eqref{eq:critical-visited-sphere-first-moment-growth}, for all sufficiently large $n$,  $\E[M_n] /2 \geq (H(R)-\varepsilon)^n$.
Therefore,
\begin{equation*}
    \pr\bigl(M_n\geq(H(R)-\varepsilon)^n\bigr)\geq c,
\end{equation*}
and hence 
\begin{equation*}
    \pr\left(\limsup_{n\to\infty}M_n^{1/n}\geq H(R)-\varepsilon\right)>0.
\end{equation*}

Changing the center of the spheres by a fixed distance does not affect $\limsup_{n\to\infty}M_n^{1/n}$. 
Hence, for every $h>0$, the event $\{\limsup_{n\to\infty}M_n^{1/n}\leq h\}$ is inherited by descendant BRWs.
The zero--one argument of \cite[Proposition~5.6]{LP17} therefore implies that $\limsup_{n\to\infty}M_n^{1/n}$ is almost surely constant on the survival event. 
Since $H(R)-\varepsilon>1$, the event $\{\limsup_{n\to\infty}M_n^{1/n}\geq H(R)-\varepsilon \}$ is contained in the survival event. 
Its positive probability therefore implies
 \begin{equation*}
    \limsup_{n\to\infty}M_n^{1/n}\geq H(R)-\varepsilon
\end{equation*}
almost surely on the survival event. 
Letting $\varepsilon\downarrow 0$ and combining with the upper bound gives \eqref{eq:M-as}. 
The logarithmic statement follows immediately from \eqref{eq:critical-visited-sphere-first-moment-growth} and \eqref{eq:M-as}.
\end{proof}

\begin{remark}
The same conclusion holds when the offspring distribution $\nu$ satisfies Assumption~\ref{ass:offspring} with mean $r\in(1,R]$. 
The case $r=R$ has been established above.
For $1<r<R$, the corresponding one-point tail estimate holds, while the subcritical triangle estimate and the Ancona inequalities yield the Green function estimates used above with $G_r$ in place of $G_R$.
Together with $H(r)>1$ and $H_n(r)\asymp H(r)^n$ from \cite[Theorem~3.1]{SWX23}, the same first- and second-moment argument yields
\begin{equation*}
    \limsup_{n\to\infty}M_n^{1/n}=H(r)
\end{equation*}
almost surely on the survival event.
\end{remark}

\section{The critical Euclidean case} \label{sec:critical-euclidean-case}
In this section, let $\Gamma=\mathbb{Z}^d$ with $d\geq5$, and let the underlying motion be the simple RW.
Suppose that the offspring distribution $\nu$ satisfies Assumption~\ref{ass:offspring} with mean $r=1$.
Since the spectral radius of the simple RW on $\mathbb{Z}^d$ is $1$, we have $R=1$. Thus $r=R$, so the BRW is critical.
Write
\begin{equation*}
    G(x,y):=G_1(x,y) \quad \text{and}\quad
    \Phi(a;x,y):=\sum_{z\in\mathbb{Z}^d}G(a,z)G(z,x)G(z,y).
\end{equation*}

For $x\in\mathbb{Z}^d$, set $\ang{x}:=1\vee |x|$, where $|x|$ denotes the $\ell^1$-norm. 
By the standard Green function estimate, for all $x,y\in\mathbb{Z}^d$,
\begin{equation} \label{eq:critical-euclidean-Green-decay}
    G(x,y)\asymp \ang{x-y}^{-d+2};
\end{equation}
see, e.g., \cite[Theorem~4.3.1]{LL10}.

\subsection{Green function estimates on $\mathbb{Z}^d$} \label{subsec:critical-euclidean-green-estimates}
We begin with the following elementary lattice sum estimates, which we will use throughout this section.

\begin{lemma} \label{lem:critical-euclidean-lattice-sum}
For every $s>0$, there exist constants $C=C(d,s)<\infty$ and $c=c(d,s)>0$ such that, for all integers $N\geq1$,
\begin{equation} \label{eq:critical-euclidean-lattice-sum}
    \sum_{z:|z|\leq N} \ang{z}^{-s} \leq C
    \begin{cases}
        N^{d-s}, & 0<s<d,\\
        1+\log N, & s=d,\\
        1, & s>d.
    \end{cases}
\end{equation}
When $0<s<d$, the matching lower bound holds:
\begin{equation} \label{eq:critical-euclidean-lattice-sum-lower}
    \sum_{z:|z|\leq N}\ang{z}^{-s} \geq c N^{d-s}.
\end{equation}
For $s>d$, we also have the tail estimate
\begin{equation} \label{eq:critical-euclidean-lattice-tail}
    \sum_{z:|z|>N}\ang{z}^{-s} \leq CN^{d-s}.
\end{equation}
\end{lemma}

\begin{proof}
For $k \geq 1$, we have $\#\{ z\in \mathbb{Z}^d :|z|=k \} \asymp k^{d-1}$.
The estimates therefore follow by comparing $1+\sum_{1\leq k\leq N}k^{d-1-s}$ and $\sum_{k>N}k^{d-1-s}$ with the corresponding integrals.
\end{proof}

We next use these lattice estimates to establish a convolution bound for power kernels.

\begin{lemma} \label{lem:critical-euclidean-power-convolution}
Let $\alpha,\beta>0$ satisfy $\alpha+\beta>d$ and $\alpha,\beta\neq d$. 
Then there exists a constant $C=C(d,\alpha,\beta)<\infty$ such that, for all $x,y\in\mathbb{Z}^d$,
\begin{equation*}
     \sum_{z\in\mathbb{Z}^d}\ang{z-x}^{-\alpha}\ang{z-y}^{-\beta} \leq C\ang{x-y}^{-\theta},
\end{equation*}
where $\theta:=\min\{\alpha,\beta,\alpha+\beta-d\}$.
\end{lemma}

\begin{proof}
The case $x=y$ follows immediately from $\alpha+\beta>d$ and Lemma~\ref{lem:critical-euclidean-lattice-sum}. 
Assume now that $x\neq y$ and set $N:=|x-y| = \ang{x-y}\geq 1$. 
Define
\begin{equation*}
    D_x:=\{z:|z-x|\leq N/2\}, \qquad D_y:=\{z:|z-y|\leq N/2\},
\end{equation*}
and let $D_{\infty}$ be the complement of $D_x\cup D_y$.

For $z\in D_x$, we have $\ang{z-y}\geq N/2$. 
Hence, by Lemma~\ref{lem:critical-euclidean-lattice-sum} and $\alpha \neq d$,
\begin{equation*}
    \sum_{z\in D_x}\ang{z-x}^{-\alpha}\ang{z-y}^{-\beta}  \leq CN^{-\beta+(d-\alpha)_+}.
\end{equation*}
Similarly, since $\beta \neq d$,
\begin{equation*}
    \sum_{z\in D_y}\ang{z-x}^{-\alpha}\ang{z-y}^{-\beta}  \leq CN^{-\alpha+(d-\beta)_+}.
\end{equation*}

For $z\in D_{\infty}$, both $\ang{z-x}$ and $\ang{z-y}$ are at least $N/2$. 
Consequently,
\begin{equation*}
  \sum_{ \substack{z\in D_\infty\\|z-x|\leq2N} } \ang{z-x}^{-\alpha}\ang{z-y}^{-\beta}  \leq CN^{d-\alpha-\beta}.
\end{equation*}
For $|z-x|>2N$, we have $\ang{z-y}\asymp\ang{z-x}$. 
Since $\alpha+\beta >d$, the tail estimate \eqref{eq:critical-euclidean-lattice-tail} yields the same upper bound. 
Therefore,
\begin{equation*}
    \sum_{z\in D_\infty}\ang{z-x}^{-\alpha}\ang{z-y}^{-\beta} \leq CN^{d-\alpha-\beta}.
\end{equation*}

Combining the three estimates and using
\begin{equation*}
    \max\{-\beta+(d-\alpha)_+,-\alpha+(d-\beta)_+,d-\alpha-\beta\} =-\theta,
\end{equation*}
we obtain the desired result.
\end{proof}

We will also use the following pointwise comparison for products of Green functions.

\begin{lemma} \label{lem:critical-euclidean-product}
There exists a constant $C=C(d)<\infty$ such that, for all $x,y,z\in\mathbb{Z}^d$,
\begin{equation*}
    G(z,x)G(z,y) \leq CG(x,y)\bigl[G(z,x)+G(z,y)\bigr].
\end{equation*}
\end{lemma}

\begin{proof}
By the triangle inequality, $\max\{\ang{z-x}, \ \ang{z-y}\}\geq \ang{x-y}/2$.
Hence, applying \eqref{eq:critical-euclidean-Green-decay}, we obtain 
\begin{equation*}
\min\{G(z,x), \ G(z,y)\} \leq C G(x,y),
\end{equation*}
which gives the desired result.
\end{proof}

We now combine the preceding estimates to obtain an upper bound for $\Phi(a;x,y)$, which in particular implies its finiteness.

\begin{corollary} \label{cor:critical-euclidean-Phi-finiteness}
There exists a constant $C=C(d)<\infty$ such that, for all $a,x,y\in\mathbb{Z}^d$,
\begin{equation} \label{eq:critical-euclidean-Phi-finiteness}
    \Phi(a;x,y)\leq CG(x,y) \left[\ang{a-x}^{-d+4}+\ang{a-y}^{-d+4} \right].
\end{equation}
In particular, $\Phi(a;x,y)<\infty$.
\end{corollary}

\begin{proof}
By Lemma~\ref{lem:critical-euclidean-product}, we have
\begin{equation*}
    \Phi(a;x,y) \leq CG(x,y)\sum_{z\in\mathbb{Z}^d}G(a,z)\bigl(G(z,x)+G(z,y)\bigr).
\end{equation*}
Applying the Green function estimate \eqref{eq:critical-euclidean-Green-decay} and Lemma~\ref{lem:critical-euclidean-power-convolution} with $\alpha=\beta=d-2$, we get
\begin{equation*}
    \sum_{z\in\mathbb{Z}^d}G(a,z)G(z,x) \leq C\ang{a-x}^{-d+4}.
\end{equation*}
The same estimate holds with $x$ replaced by $y$. 
Combining these estimates, we obtain the desired bound, which also shows that $\Phi(a;x,y)<\infty$ for all $a,x,y\in\mathbb{Z}^d$.
\end{proof}

We next sharpen this upper bound to determine the precise order of $\Phi(a;x,y)$.

\begin{proposition} \label{prop:critical-euclidean-Phi-asymptotic}
There exist constants $c=c(d)>0$ and $C=C(d)<\infty$ such that, for all $a, x, y\in \mathbb{Z}^d$, 
\begin{equation} \label{eq:critical-euclidean-Phi-asymptotic}
    c l_1^{-d+4} l_2^{-d+2} \leq \Phi(a;x,y) \leq C l_1^{-d+4} l_2^{-d+2},
\end{equation}
where $l_1\leq l_2\leq l_3$ is the increasing rearrangement of $\ang{x-a}$, $ \ang{y-a}$ and $\ang{x-y}$.
\end{proposition}

\begin{proof}
Let $(u,v,w)$ be a permutation of $(a,x,y)$ such that $\ang{u-v}=l_1$.
The other two mutual distances are $l_2$ and $l_3$. 
By the triangle inequality,
\begin{equation} \label{eq:critical-euclidean-distance-comparison}
    \ang{u-w}\asymp\ang{v-w}\asymp l_2.
\end{equation}
Moreover, by the symmetry of $G$, we have
\begin{equation} \label{eq:Phi-uvw}
	\Phi(a;x,y) =\sum_{z\in\mathbb{Z}^{d}}G(z,u)G(z,v)G(z,w) =\Phi(v;u,w).
\end{equation}

For the upper bound, applying Corollary~\ref{cor:critical-euclidean-Phi-finiteness} to \eqref{eq:Phi-uvw} and using \eqref{eq:critical-euclidean-distance-comparison}, we obtain
\begin{align*}
    \Phi(a;x,y) &\leq CG(u,w)\left(\ang{v-u}^{-d+4} + \ang{v-w}^{-d+4} \right) \\
    &\leq C l_2^{-d+2}\left( l_1^{-d+4}+l_2^{-d+4}\right) \\
    &\leq C l_1^{-d+4} l_2^{-d+2}.
\end{align*}

For the lower bound, first suppose that $l_1\leq 4$.  
Taking only the term $z=u$ in \eqref{eq:Phi-uvw} and using \eqref{eq:critical-euclidean-Green-decay}, we obtain
\begin{equation*}
    \Phi(a;x,y) \geq c l_1^{-d+2} l_2^{-d+2} \geq c l_1^{-d+4} l_2^{-d+2},
\end{equation*}
where the last inequality follows from  $l_1\leq 4$.

Suppose now that $l_1>4$ and consider 
\begin{equation*}
    B_u:=\{z\in\mathbb{Z}^d:|z-u|\leq l_1/4\}.
\end{equation*}
For $z\in B_u$,  the triangle inequality and \eqref{eq:critical-euclidean-distance-comparison} imply $ \ang{z-v}\asymp l_1$ and $ \ang{z-w}\asymp l_2$. 
Hence, by \eqref{eq:Phi-uvw} and  \eqref{eq:critical-euclidean-Green-decay}, we obtain
\begin{equation*}
    \Phi(a;x,y)  \geq c l_1^{-d+2} l_2^{-d+2}  \sum_{z\in B_u}\ang{z-u}^{-d+2} \geq c l_1^{-d+4} l_2^{-d+2},
\end{equation*}
where we applied \eqref{eq:critical-euclidean-lattice-sum-lower} with $s=d-2$ in the last inequality.
This completes the proof.	
\end{proof}

\subsection{Verification of the Green function criterion}

We first verify the triangle condition~\eqref{eq:criterion-triangle}.

\begin{lemma} \label{lem:critical-euclidean-triangle}
There exists a constant $C=C(d)<\infty$ such that, for all $a, x \in\mathbb{Z}^d$,
\begin{equation} \label{eq:critical-euclidean-triangle}
    \sum_{z\in\mathbb{Z}^d} G(a,z)^2 G(z,x) \leq CG(a,x).
\end{equation}
\end{lemma}

\begin{proof}
By \eqref{eq:critical-euclidean-Green-decay}, we have
\begin{equation*}
    \sum_{z\in\mathbb{Z}^d}G(a,z)^2 G(z,x) \leq C\sum_{z\in\mathbb{Z}^d}\ang{z-a}^{-2d+4}\ang{z-x}^{-d+2}.
\end{equation*}
Since $2d-4>d$, applying Lemma~\ref{lem:critical-euclidean-power-convolution} with $\alpha=2d-4$ and $\beta=d-2$, we obtain
\begin{equation*}
    \sum_{z\in\mathbb{Z}^d}G(a,z)^2G(z,x) \leq C\ang{a-x}^{-d+2} \leq CG(a,x).
\end{equation*}
This completes the proof.
\end{proof}

\begin{remark}  \label{rem:critical-euclidean-dimension-threshold}
The restriction $d\geq5$ ensures that $2(d-2)>d$, or equivalently, that the Green function is square summable.
We use this property in the proofs of Corollary~\ref{cor:critical-euclidean-Phi-finiteness} and Lemma~\ref{lem:critical-euclidean-triangle}.

In dimension $4$, the Green function estimate~\eqref{eq:critical-euclidean-Green-decay} remains valid. 
For $|x|\geq2$, we restrict the sum to $|z|\leq |x|/2$ and obtain
\begin{equation*}
     \sum_{z\in\mathbb{Z}^4}G(0,z)^2 G(z,x)  \geq c G(0,x) \log |x|.
\end{equation*}
Thus, the triangle condition~\eqref{eq:criterion-triangle} fails in dimension $4$.
\end{remark}

We now turn to the mixed convolution condition~\eqref{eq:criterion-mixed}.

\begin{proposition} \label{prop:critical-euclidean-G-Phi-convolution}
There exists a constant $C=C(d)<\infty$ such that, for all $a,x,y\in\mathbb{Z}^d$,
\begin{equation} \label{eq:critical-euclidean-G-Phi-convolution}
    \sum_{z\in\mathbb{Z}^d}G(a,z)G(z,x)\Phi(z;x,y)  \leq C\Phi(a;x,y).
\end{equation}
\end{proposition}

\begin{proof}
Expanding $\Phi$ and applying Tonelli's theorem, we have 
\begin{equation} \label{eq:expanding-sum}
    \sum_{z\in\mathbb{Z}^d}G(a,z)G(z,x)\Phi(z;x,y)= \sum_{z,u\in\mathbb{Z}^d} G(a,z)G(z,x)G(z,u)G(u,x)G(u,y).
\end{equation}
By the symmetry of $G$ and Lemma~\ref{lem:critical-euclidean-product}, we get
\begin{equation*}
	G(z,x)G(u,x) \leq CG(z,u)\bigl(G(z,x)+G(u,x)\bigr).
\end{equation*}
Hence, the sum in \eqref{eq:expanding-sum} is bounded by $C(U_1+U_2)$, where
\begin{equation*}
    U_1:=\sum_{z,u\in\mathbb{Z}^d}G(a,z)G(z,x)G(z,u)^2G(u,y)
\end{equation*}
and
\begin{equation*}
    U_2:=\sum_{z,u\in\mathbb{Z}^d}G(a,z)G(z,u)^2G(u,x)G(u,y).
\end{equation*}

Applying Lemma~\ref{lem:critical-euclidean-triangle} to the sum over $u$, we obtain
\begin{equation*}
    U_1\leq C\sum_{z\in\mathbb{Z}^d} G(a,z)G(z,x)G(z,y)=C\Phi(a;x,y).
\end{equation*}
Similarly, applying Lemma~\ref{lem:critical-euclidean-triangle} to the sum over $z$ and using the symmetry of $G$, we obtain
\begin{equation*}
    U_2 \leq C\sum_{u\in\mathbb{Z}^d} G(a,u)G(u,x)G(u,y) =C\Phi(a;x,y).
\end{equation*}
Combining these estimates gives \eqref{eq:critical-euclidean-G-Phi-convolution}.
\end{proof}

\begin{proof}[Proof of Corollary~\ref{cor:critical-euclidean-two-point-tail-high}]
The Green function is finite by \eqref{eq:critical-euclidean-Green-decay}. 
Lemma~\ref{lem:critical-euclidean-triangle} and Proposition~\ref{prop:critical-euclidean-G-Phi-convolution} verify the conditions~\eqref{eq:criterion-triangle} and~\eqref{eq:criterion-mixed}, respectively.
Theorem~\ref{thm:two-point-tail-estimate} therefore yields the two-point tail estimate.
Finally, Proposition~\ref{prop:critical-euclidean-Phi-asymptotic} gives the stated bounds in terms of $l_1$ and $l_2$.
This completes the proof.
\end{proof}

\subsection{An application to visiting probabilities of finite sets}
We conclude this section by applying the mixed moment bound to visiting probabilities of finite sets.
For a finite nonempty set $A\subset\mathbb{Z}^d$, write 
\begin{equation*}
	G(a, A):=\sum_{x\in A} G(a,x) \quad \text{and} \quad  N_A:=\sum_{x\in A}L(x).
\end{equation*}
Then $\{N_A>0\}$ is precisely the event that the BRW visits $A$.

We first establish a general lower bound.

\begin{proposition} \label{prop:critical-euclidean-finite-set-visiting}
There exists a constant $c=c(d,\nu)>0$ such that, for all $a\in\mathbb{Z}^d$ and all finite nonempty sets $A\subset\mathbb{Z}^d$,
\begin{equation*} 
     \pr_a\bigl(N_A>0\bigr) \geq c \frac{G(a,A)^2}{\sum_{z\in\mathbb{Z}^d}G(a,z)G(z,A)^2}.
\end{equation*}
\end{proposition}

\begin{proof}
By the many-to-one formula, we have $\E_a[N_A] =G(a,A)$. 
Moreover, Corollary~\ref{cor:mixed-moment-upper-bound} with $p=q=1$ gives
\begin{align*}
    \E_a[N_A^2] &=\sum_{x,y\in A}\E_a[L(x)L(y)] \\
    &\leq C\sum_{x,y\in A}\Phi(a;x,y)=C\sum_{z\in\mathbb{Z}^d}G(a,z)
    \left(\sum_{x\in A}G(z,x)\right)^2 =C\sum_{z\in\mathbb{Z}^d}G(a,z)G(z,A)^2.
\end{align*}
Therefore, applying the Paley--Zygmund inequality, we obtain
\begin{equation*}
	  \pr_a(N_A>0) \geq \frac{\E_a[N_A]^2}{\E_a[N_A^2]} \geq c \frac{G(a,A)^2}{\sum_{z\in\mathbb{Z}^d}G(a,z)G(z,A)^2}.
\end{equation*}
This completes the proof.
\end{proof}

We next apply Proposition~\ref{prop:critical-euclidean-finite-set-visiting} to balls contained in coordinate subspaces.
For $1\leq k\leq d$ and $M\geq2$, let
\begin{equation*}
    B_k(M):=\bigl\{(u,0)\in\mathbb{Z}^k\times\mathbb{Z}^{d-k}:\|u\|\leq M\bigr\},
\end{equation*}
where $\| \cdot \|$ denotes the Euclidean norm. 
For $a\notin B_k(M)$, set $s:=\inf_{x\in B_k(M)}\|a-x\|$.

We first consider $s\geq M$ and $1\leq k\leq d-3$.
By the Green function estimate~\eqref{eq:critical-euclidean-Green-decay} and the equivalence of the Euclidean and $\ell^1$ norms, 
\begin{equation} \label{eq:critical-euclidean-flat-ball-first-moment}
    G \left(a,B_k(M) \right) \asymp  s^{-d+2} M^k .
\end{equation}
We also have 
\begin{equation} \label{eq:critical-euclidean-flat-ball-second-moment-kernel}
    \sum_{z\in\mathbb{Z}^d}G(a,z) G\left(z,B_k(M)\right)^2 \leq C s^{-d+2}
    \begin{cases}
        M^{k+1}, & k=d-3,\\
        M^k\log M, & k=d-4,\\
        M^k, & 1\leq k\leq d-5.
    \end{cases}
\end{equation}
To see this, note that for $x,y\in B_k(M)$ with $s\geq M$, we have $ \ang{a-x}\asymp\ang{a-y}\asymp s$ and  $\ang{x-y}\leq CM$.
Hence, by translation invariance and Proposition~\ref{prop:critical-euclidean-Phi-asymptotic},
\begin{equation*}
    \Phi(a;x,y)=\Phi(0;x-a,y-a)\leq Cs^{-d+2}\ang{x-y}^{-d+4}.
\end{equation*}
Summing over $x, y\in B_k(M)$, we obtain
\begin{align} 
    &\sum_{z\in\mathbb{Z}^d}G(a,z)G\bigl(z,B_k(M)\bigr)^2 =\sum_{x,y\in B_k(M)}\Phi(a;x,y)  \notag\\
    &\qquad \leq Cs^{-d+2}\sum_{x,y\in B_k(M)}\ang{x-y}^{-d+4}
    \leq Cs^{-d+2}M^k  \sum_{\substack{u\in\mathbb{Z}^k\\\|u\|\leq2M}} \ang{u}^{-d+4}. \label{eq:1}
\end{align}
By the same lattice sum argument as in Lemma~\ref{lem:critical-euclidean-lattice-sum}, we have 
\begin{equation*}
    \sum_{\substack{u\in\mathbb{Z}^k\\\|u\| \leq2M}}
    \ang{u}^{-d+4} \leq C
    \begin{cases}
        M, & k=d-3, \\
        \log M, & k=d-4,\\
        1, & 1\leq k\leq d-5.
    \end{cases}
\end{equation*}
Thus, combining this with \eqref{eq:1} gives the desired estimate \eqref{eq:critical-euclidean-flat-ball-second-moment-kernel}.

The estimates in \eqref{eq:critical-euclidean-flat-ball-first-moment} and \eqref{eq:critical-euclidean-flat-ball-second-moment-kernel}, together with Proposition~\ref{prop:critical-euclidean-finite-set-visiting} and a suitable subball argument, yield the following corollary.

\begin{corollary} \label{cor:critical-euclidean-flat-ball-visiting}
There exists a constant $c=c(d,\nu)>0$ such that, for all $1\leq k\leq d$, $M\geq2$, and all $a\notin B_k(M)$ with $ s:=\inf_{x\in B_k(M)}\|a-x\| \geq2$, 
\begin{equation}\label{eq:critical-euclidean-flat-ball-visiting}
    \pr_a\bigl(N_{B_k(M)}>0\bigr)
    \geq cs^{-d+2}
    \begin{cases}
        t^{d-4}, & k\geq d-3,\\
        t^{d-4}(\log t)^{-1}, & k=d-4,\\
        t^k, & 1\leq k\leq d-5,
    \end{cases}
\end{equation}	
where $t:=\min\{s,M\}$.
\end{corollary}

\begin{proof}
First suppose that $s\geq M$ and $k\leq d-3$.
Then $t=M$, and the claimed bounds follow from Proposition~\ref{prop:critical-euclidean-finite-set-visiting}, \eqref{eq:critical-euclidean-flat-ball-first-moment}, and \eqref{eq:critical-euclidean-flat-ball-second-moment-kernel}.

To handle the remaining cases, set $j:=\min \{k, d-3\}$ and choose $b \in B_k(M)$ such that $\|a-b\|=s$.  
Let $t_0:=\max\{ 16, 4\sqrt{d} \}$.

Suppose that $t\leq t_0$. Taking $A=\{b\}$ in Proposition~\ref{prop:critical-euclidean-finite-set-visiting} and using \eqref{eq:critical-euclidean-triangle}, we obtain
\begin{equation*}
	 \pr_a(N_{B_k(M)}>0) \geq \pr_a(L(b)>0) \geq cG(a,b) \geq cs^{-d+2}.
\end{equation*}
Since $2\leq t\leq t_0$, this gives the claimed bounds after adjusting the constant $c$.

Now suppose that $t>t_0$. 
Choose $b_0 \in\mathbb{Z}^k \times \{0\}^{d-k}$ such that
\begin{equation*}
    \left\|b_0-\left(1-\frac{t}{2M}\right)b\right\| \leq \frac{\sqrt d}{2}\leq\frac{t}{8}.
\end{equation*}
Then $ \|b_0\| \leq M-3t/8$ and $ \|b_0-b\| \leq 5t/8$.
Consequently, the translated coordinate ball $B':=b_0+B_j(t/8)$ satisfies
\begin{equation*}
    B'\subseteq B_k(M) \qquad \text{and} \qquad s\leq s':= \inf_{x\in B'}\|a-x\|\leq s+\frac{5t}{8}\leq 2s.
\end{equation*}
Since $2\leq t/8 \leq s'$ and $j \leq d-3$, we can apply the bounds established in the first case to $B'$ by translation invariance.
Using $B'\subseteq B_k(M)$, $s'\asymp s$, and $\log(t/8)\asymp\log t$, we obtain the claimed bounds.
\end{proof}

Under our offspring assumptions,  Corollary~\ref{cor:critical-euclidean-flat-ball-visiting} recovers the lower bounds in \cite[Theorem~21]{Zhu26}, while the second-moment estimate above corresponds to \cite[Lemma~22]{Zhu26}.

\section*{Acknowledgement}
The authors would like to thank Professor  Wei Xu for his valuable advice and guidance, which greatly improved the quality of this paper, and Beijing Institute of Technology for providing a stimulating academic environment.

\bibliographystyle{plainnat}
\bibliography{tail}

@book {Woe00,
    AUTHOR = {Woess, Wolfgang},
     TITLE = {Random walks on infinite graphs and groups},
    SERIES = {Cambridge Tracts in Mathematics},
    VOLUME = {138},
 PUBLISHER = {Cambridge University Press, Cambridge},
      YEAR = {2000},
     PAGES = {xii+334},
      ISBN = {0-521-55292-3},
   MRCLASS = {60B15 (60G50 60J10)},
  MRNUMBER = {1743100},
MRREVIEWER = {Donald I. Cartwright},
       URL = {https://doi.org/10.1017/CBO9780511470967},
}

@article{AHJ21,
  title={On the tail of the branching random walk local time},
  author={Angel, Omer and Hutchcroft, Tom and J{\'a}rai, Antal},
  journal={Probability Theory and Related Fields},
  volume={180},
  number={1},
  pages={467--494},
  year={2021},
  publisher={Springer}
}

@article{BP94,
  title={Markov chains indexed by trees},
  author={Benjamini, Itai and Peres, Yuval},
  journal={The Annals of Probability},
  pages={219--243},
  year={1994},
  publisher={JSTOR}
}

@article{GM06, author = {Gantert, Nina and M{\"u}ller, Sebastian}, 
title = {The critical branching {M}arkov chain is transient}, 
journal = {Markov Processes and Related Fields}, 
volume = {12}, 
pages = {805--814}, 
year = {2006} }

@book{LP17,
  title={Probability on trees and networks},
  author={Lyons, Russell and Peres, Yuval},
  volume={42},
  year={2017},
  publisher={Cambridge University Press}
}

@article{Gou15,
  title={Martin boundary of random walks with unbounded jumps in hyperbolic groups},
  author={Gou{\"e}zel, S{\'e}bastien},
  journal={The Annals of Probability},
  pages={2374--2404},
  year={2015},
  publisher={JSTOR}
}

@book{CDP90,
  author    = {Coornaert, Michel and Delzant, Thomas and Papadopoulos, Athanase},
  title     = {G{\'e}om{\'e}trie et th{\'e}orie des groupes: les groupes hyperboliques de Gromov},
  series    = {Lecture Notes in Mathematics},
  volume    = {1441},
  publisher = {Springer-Verlag},
  address   = {Berlin},
  year      = {1990},
  doi       = {10.1007/BFb0084913}
}

@book{LL10,
  author    = {Lawler, Gregory F. and Limic, Vlada},
  title     = {Random Walk: A Modern Introduction},
  series    = {Cambridge Studies in Advanced Mathematics},
  volume    = {123},
  publisher = {Cambridge University Press},
  address   = {Cambridge},
  year      = {2010}
}

@article{SWX23,
  author  = {Sidoravicius, Vladas and Wang, Longmin and Xiang, Kainan},
  title   = {Limit set of branching random walks on hyperbolic groups},
  journal = {Communications on Pure and Applied Mathematics},
  volume  = {76},
  number  = {10},
  pages   = {2765--2803},
  year    = {2023},
  doi     = {10.1002/cpa.22088}
}

@article{GL26,
  author  = {Gao, Ruoxiang and Lai, Shuwen},
  title   = {On the Tail of the Local Time of Branching Random Walks on Finitely Generated Groups},
  journal = {Electronic Communications in Probability},
  year    = {2026},
  note    = {To appear}
}

@article{Zhu26,
  author  = {Zhu, Qingsan},
  title   = {On the Critical Branching Random Walk {I}: Branching Capacity and Visiting Probability},
  journal = {Probability Theory and Related Fields},
  year    = {2026},
  doi     = {10.1007/s00440-026-01530-9}
}

@article{CL26,
  author  = {Chen, Qitong and Lai, Shuwen},
  title   = {Sharp Two-Point Local Time Estimates for Critical Branching Random Walks in Low Dimensions},
  journal = {Statistics \& Probability Letters},
  volume  = {239},
  pages   = {110888},
  year    = {2026},
  doi     = {10.1016/j.spl.2026.110888}
}

@article{ASS26,
  author  = {Asselah, Amine and Schapira, Bruno and Sousi, Perla},
  title   = {Local times and capacity for transient branching random walks},
  journal = {Probability Theory and Related Fields},
  volume  = {194},
  pages   = {1443--1493},
  year    = {2026},
  doi     = {10.1007/s00440-025-01378-5}
}

@article{Zhu17,
  author  = {Zhu, Qingsan},
  title   = {On the critical branching random walk {II}: Branching capacity and branching recurrence},
  journal = {arXiv preprint arXiv:1612.00161},
  year    = {2017}
}

@article{BZ07,
  author  = {Birkner, Matthias and Z{\"a}hle, Iljana},
  title   = {A functional {CLT} for the occupation time of a state-dependent branching random walk},
  journal = {The Annals of Probability},
  volume  = {35},
  number  = {6},
  pages   = {2063--2090},
  year    = {2007},
  doi     = {10.1214/009117907000000150}
}

@article{LZ10,
  author  = {Lalley, Steven P. and Zheng, Xinghua},
  title   = {Spatial Epidemics and Local Times for Critical Branching Random Walks in Dimensions 2 and 3},
  journal = {Probability Theory and Related Fields},
  volume  = {148},
  number  = {3--4},
  pages   = {527--566},
  year    = {2010}
}

@article{LZ11,
  author  = {Lalley, Steven P. and Zheng, Xinghua},
  title   = {Occupation Statistics of Critical Branching Random Walks in Two or Higher Dimensions},
  journal = {The Annals of Probability},
  volume  = {39},
  number  = {1},
  pages   = {327--368},
  year    = {2011}
}

@article{LGL15,
  author  = {Le Gall, Jean-Fran{\c{c}}ois and Lin, Shen},
  title   = {The Range of Tree-Indexed Random Walk in Low Dimensions},
  journal = {The Annals of Probability},
  volume  = {43},
  number  = {5},
  pages   = {2701--2728},
  year    = {2015}
}

@unpublished{Lai26,
  author = {Lai, Shuwen},
  title  = {Two-sided estimates for two-point local time tails of critical branching random walks in dimension four},
  year   = {2026},
  note   = {Manuscript}
}
\end{document}